\documentstyle[12pt,amscd]{amsart}

\theoremstyle{plain} 

\newtheorem{ThNum}{Theorem}

\newtheorem{LemNum}{Lemma}
\newtheorem{Remark}{Remark}

\newcommand{\hyp}{{\bold H}}

\newcommand{\real}{{\bold R}}
\newcommand{\complex}{{\bold C}}
\newcommand{\projective}{{\bold P}}

\newcommand{\PGL}{\operatorname{PGL}} 
 
\newcommand{\SL}{\operatorname{SL}}

\newcommand{\trace}{\operatorname{trace}} 
 
\newcommand{\Area}{\operatorname{Area}}

\begin{document}

\title[Configuration spaces and hyperbolic Dehn fillings]{%
        Configuration spaces of points on the circle and 
	hyperbolic Dehn fillings} 

\author[S. Kojima]{%
        Sadayoshi Kojima} 
\address{%
        Department of Mathematical and Computing Sciences \\
        Tokyo Institute of Technology \\
        Ohokayama, Meguro \\    
        Tokyo 152-8552 Japan} 
\email{%
        sadayosi@@is.titech.ac.jp}  
\author[H. Nishi]{%
        Haruko Nishi} 
\address{%
        Department of Mathematics \\
        Kyushu University \\ 
        33, Fukuoka 812-8581 Japan} 
\email{%
        nishi@@math.kyushu-u.ac.jp} 
\author[Y. Yamashita]{%
        Yasushi Yamashita} 
\address{%
        Department of Information and Computer Sciences \\
        Nara Women's University \\ 
        Kita-Uoya Nishimachi \\ 
        Nara 630-8506, Japan} 
\email{%
        yamasita@@ics.nara-wu.ac.jp} 




\begin{abstract} 
A purely combinatorial compactification of 
the configuration space of  $n \, (\geq 5)$  distinct points 
with equal weights 
in the real projective line was introduced by M. Yoshida.  
We geometrize it so that it will be a real 
hyperbolic cone-manifold of finite volume with dimension  $n - 3$.   
Then, we vary weights for points.  
The geometrization still makes sense and yields a deformation. 
The effectivity of deformations arisen in this manner 
will be locally described 
in the existing deformation theory of hyperbolic structures 
when  $n - 3 = 2, 3$.  
\end{abstract} 

\maketitle


\section{Introduction}

Let  $X(n)$  be a space of configurations of $n$ distinct 
points in the real projective line $\real \projective^1$  up to 
projective automorphisms.  
$X(n)$  can be expressed by the point set,   
\begin{equation*} 
	X(n) = ((\real \projective^1)^n - {\bold D}) 
		/ \PGL(2,\real),  
\end{equation*} 
where  ${\bold D}$  is the big diagonal set, 
\begin{equation*} 
	{\bold D} = \{ (\alpha_1, \cdots, \alpha_n) \in 
		(\real \projective)^n \, \vert \, 
		 \alpha_i = \alpha_j \; \; \text{for some  $i \ne j$} \},    
\end{equation*}  
and  $\PGL(2, \real)$  acts on  $(\real \projective^1)^n$  
diagonally.  
We assume that the number of points is at 
least five throughout this paper 
(for the case n=4, see \cite{MorinNishi}). 

There are two obvious observations.   
$X(n)$  is not connected since we are not allowed 
to have collisions of points.  
Each component can be labeled by 
a circular permutation of $n$ letters up to reflection.  
In particular, the number of connected components is  $(n-1)!/2$. 
Also a configuration can be normalized by sending 
three consecutive points to  $\{ 0, 1, \infty \}$  
so that the other points lie in the open unit interval  $(0, 1)$.  
Hence each component of  $X(n)$  can be 
identified with the set of ordered  $n-3$  points in  $(0, 1)$, 
and in particular is 
homeomorphic to a cell of dimension  $n-3$.  

M. Yoshida introduced  
a purely combinatorial compactification  
$X_Y(n)$  of  $X(n)$ 
in \cite{Yoshida} by deleting the following smaller set  
${\bold D}_*$  instead of the big diagonal  ${\bold D}$,  
\begin{equation*} 
	{\bold D}_* = 
	\{ (\alpha_1, \cdots, \alpha_n) 
		\in ({\real \projective}^1)^n \, \vert \, 
	\alpha_{i_1} = \alpha_{i_2} = \cdots = \alpha_{i_{[\frac{n+1}2]}} 
		\; \;  
	\text{for some distict $i_1, \cdots, i_{[\frac{n+1}2]}$}\}
\end{equation*}
where  $[x]$  denotes the maximal integer which does not 
exceed  $x$.  
When  $n$  is even,  
$X_Y(n)$  is not compact in fact,  
but there is a natural combinatorial interpretation of the ends.  

On the other hand, 
Thurston gave a family of incomplete complex hyperbolic 
metrics on the space of configurations of  $n$ distinct points 
on the complex projective line $\complex \projective^1$  in  
\cite{Thurston}.   
Their completions become complex hyperbolic 
cone-manifolds of complex dimension  $n-3$.  

By considering the real slice 
of Thurston's complex hyperbolization, 
we can give a natural real hyperbolic polyhedral structure on 
each component of  $X(n)$.    
Moreover the boundary of each polyhedron will be coded by 
degenerate configurations so that the hyperbolization 
gives rise to a geometric 
identification of pairs of faces on the components.  
The first purpose of this paper is to show 
that such identifications make  $X(n)$  lie 
as an open dense subset in    
a connected real hyperbolic cone-manifold homeomorphic to 
Yoshida's compactification in Theorem 1.  

We will see more precisely that 
the resultant of face identification is compact if  $n$  is odd, and 
noncompact of finite volume otherwise.  
If  $ n \leq 6$, it is nonsingular. 
If  $n \geq 7$,  it is singular and 
its cone angle  $\theta_n$  along a codimension two 
singular stratum satisfies the identity  
$\cos (\theta_n/6) = 1/(2 \cos (2\pi/n))$.  
Hence  $\theta_n$  lies in 
$(\pi, 2\pi)$  and approaches  $2\pi$  
as  $n \to \infty$.  

Some related works on geometric 
interpretation of  $X(n)$  from 
more analytic viewpoints can be found in 
Ap\'ery and Yoshida \cite{AperyYoshida}, 
Yamazaki and Yoshida \cite{YamazakiYoshida} 
and Yoshida \cite{YoshidaPreprint}.  

The second purpose of this paper is 
to relate our geometrization with the existing 
deformation theory of hyperbolic structures 
by relaxing the construction when  $n = 5, 6$.  
The discussion for the first purpose is based on 
the hypothesis that the points on a configuration all have 
equal weights.  
We will see what happens if we perturb weights of points slightly.  

When  $n =5$,  the geometrized configuration space is 
a nonorientable hyperbolic surface 
homeomorphic to a connected sum of five copies of 
the projective plane,  ${\#}^5 \real \projective^2$.  
We will assign a deformed hyperbolic structure  
on  ${\#}^5 \real \projective^2$  to each perturbation of weights. 
The freedom for varying 
weights on five points is of 4 dimensional, 
and the space of marked hyperbolic structures 
on  ${\#}^5 \real \projective^2$  is 
homeomorphic to  $\real^9$.   
We will show 
that this assignment is differentiable 
and that the derivative at equal weights has full rank  
in Theorem 2.  

When  $n = 6$,  the geometrized configuration space  
is a complete hyperbolic 3-manifold 
of finite volume with ten cusps.  
A small perturbation gives rise to not a deformation 
in the usual sense, 
but a resultant of hyperbolic Dehn filling.  
The space of hyperbolic Dehn fillings of a complete 
hyperbolic 3-manifold can be locally identified with 
the space of representations of a fundamental group up to 
conjugacy.  
It has a structure of an algebraic variety of 
complex dimension $= \text{the number of cusps}$, 
and is smooth at the complete structure.  
Hence in our case,  the space of Dehn fillings 
is locally biholomorphic to  $\complex^{10}$.    
Then we will show again that the assignment is differentiable 
and the derivative at equal weights has full rank  
in Theorem 3. 

We carry out the geometrization in the next section, 
and relate it with the deformation theory 
in the section after. 

The first author would like to thank Junjiro Noguchi 
for his comment in early stage of this work, and 
the authors all would like to thank Masaaki Yoshida for 
motivating us and for interests to this work.


\section{Geometrization}

\subsection{Hyperbolic polyhedral structure on $X(n)$}
\label{subsec:PolyStr}

To geometrize  $X(n)$, we consider Euclidean
$n$-gons with vertices marked by integers from $1$ to $n$,  
where the marking may not be cyclically monotone.  
Let $X_{n,c}$ be the set of all marked
equiangular $n$-gons up to mark preserving 
(possibly orientation reversing) congruence,  
and $X_n$ a further quotient of  $X_{n,c}$  by similarities.

For any $\alpha \in (\real \projective^1)^n - {\bold D}$, 
we assign the unit
disc in $\complex$ with $n$ points specified on the boundary.  
By the Schwarz--Christoffel mapping or its complex conjugate, 
we can map $\alpha$ to an 
{\it equiangular} $n$-gon up to mark preserving
similarity --- i.e. an element of $X_n$.  
This induces a map 
from $X(n)$ to $X_n$ 
since a projective transformation on the unit disc 
does not change the image of the map.  
It is also injective because if two
configurations $\alpha$ and $\beta$ map to the same
element of $X_{n,c}$ by $f_{\alpha}$ and $f_{\beta}$  
then $f_{\beta}\circ f_{\alpha}^{-1}$ is a
mark preserving projective automorphism of the unit disk.  
By the Carath\'{e}odory theorem, this map is surjective.  
Therefore we have proved 

\begin{LemNum}\label{Lem:Bijective}
There is a canonical homeomorphism between  $X(n)$  and  $X_n$.
\end{LemNum}

The aim of this subsection is to construct a hyperbolic polyhedral
structure on a component of $X(n)$ through the above
identification by taking the real part of Thurston's geometrization 
in \cite{Thurston}.  
Similar discussions with 
the rest of this subsection can be found also in 
\cite{BavardGhys, KojimaYamashita}. 

As we have noted in the introduction, the connected component of
$X(n)$ consists of the configurations with 
a fixed circular permutation of markings up to reflection.  
For the simplicity of index, 
we shall concentrate on the component $U$ of $X_{n,c}$
which is labeled by  $12 \cdots n$.  
Let  $U_s \subset X_n$ be the set of its mark preserving similarity
classes.  
Note that $U_s$ corresponds to the component of $X(n)$  
labeled by  $12 \cdots n$  also.
We identity $U_s$ with the set $U_1$ which, by definition, consists of
the set of equiangular polygons having its area one.  We describe
polygons as follows.

The elements of $U$ can be described by the vector of side lengths
$(x_1, x_2,\ldots, x_n)$ where $x_j$ is the length of the edge between
the vertices marked by  $j$ and $j+1$.  
Since they represent an equiangular $n$-gon, 
they satisfy:
$$
 x_1 + x_2 \zeta_n + \cdots + x_n \zeta_n^{n-1} = 0
$$
where $\zeta_n = \exp(2\pi i/n)$. We set
\begin{align*}
  {\cal E}_n & :=  \left\{(x_1, x_2,\ldots, x_n)\mid
     x_1 + x_2 \zeta_n + \cdots + x_n \zeta_n^{n-1} = 0\right\}, \\
  {\cal E}_n^+ & :=  {\cal E}_n \cap \bigcap_{j=1}^n \{x_j >0\}.
\end{align*}
Note that $U = {\cal E}_n^+$.
\begin{figure}[ht]
\begin{center}
\setlength{\unitlength}{0.00083300in}%
\begingroup\makeatletter\ifx\SetFigFont\undefined%
\gdef\SetFigFont#1#2#3#4#5{%
  \reset@font\fontsize{#1}{#2pt}%
  \fontfamily{#3}\fontseries{#4}\fontshape{#5}%
  \selectfont}%
\fi\endgroup%
\begin{picture}(2797,2905)(879,-2183)
\put(1426,614){Im}
\thicklines
\put(901,-1561){\vector( 1, 0){2700}}
\put(1501,-2161){\vector( 0, 1){2700}}
\put(1726,-1711){$x_8$}
\put(2626,-1186){$x_1$}
\put(3076,-586){$x_2$}
\put(2776,-61){$x_3$}
\put(1951,314){$x_4$}
\put(1201,-61){$x_5$}
\put(901,-811){$x_6$}
\put(1051,-1411){$x_7$}
\put(1576,-1486){$i_8$}
\put(1951,-1486){$i_1$}
\put(2776,-661){$i_2$}
\put(2776,-436){$i_3$}
\put(2326, 14){$i_4$}
\put(1801, 14){$i_5$}
\put(1276,-436){$i_6$}
\put(1276,-1261){$i_7$}
\put(3676,-1636){Re}
\put(1501,-1561){\line( 1, 0){600}}
\put(2101,-1561){\line( 1, 1){900}}
\put(3001,-661){\line( 0, 1){300}}
\put(3001,-361){\line(-1, 1){600}}
\put(2401,239){\line(-1, 0){600}}
\put(1801,239){\line(-1,-1){600}}
\put(1201,-361){\line( 0,-1){900}}
\put(1201,-1261){\line( 1,-1){300}}
\end{picture}
\end{center}
\caption{An equiangular octagon} \label{Fig:Equiangle}
\end{figure}
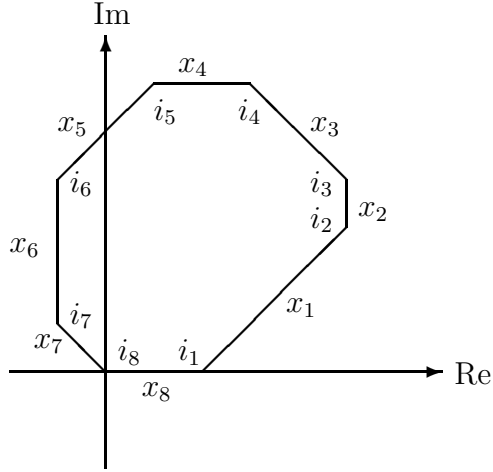

For each element $P$ of ${\cal E}_n^+$, we denote by $\Area(P)$
the area of $P$.  It is a function from ${\cal E}_n^+$ to $\real$.

\begin{LemNum}\label{Lem:Quadratic}
The function $\text{\em Area}$ is extended to a quadratic 
form of signature $(1,n-3)$ on ${\cal E}_n$.
\end{LemNum}

\begin{pf}
Suppose that $P$ is an element of $U$.  Place the edge $x_n$ of $P$ on
the real axis of complex plane and extend other edges until they touch
the real axis.  We set edge length $y_j, z_k, a_l$ and triangles
$B_s, C_t$ as in the figure. (If n is odd, then no $A_{n/2}$ appears.)
Let $A$ be the triangle $P \cup_s B_s \cup_t C_t$.

\begin{figure}[ht]
\begin{center}
\setlength{\unitlength}{0.00083300in}%
\begingroup\makeatletter\ifx\SetFigFont\undefined%
\gdef\SetFigFont#1#2#3#4#5{%
  \reset@font\fontsize{#1}{#2pt}%
  \fontfamily{#3}\fontseries{#4}\fontshape{#5}%
  \selectfont}%
\fi\endgroup%
\begin{picture}(4694,2744)(129,-2483)
\put(976,-1636){$C_2$}
\thicklines
\put(1801,-2461){\vector( 0, 1){2700}}
\put(151,-1861){\line( 1, 0){4650}}
\put(301,-1861){\line( 1, 1){2100}}
\put(2401,239){\line( 1,-1){2100}}
\put(3301,-961){\line( 0,-1){900}}
\put(1501,-1561){\line( 0,-1){300}}
\put(1340,-2154){\vector( 2, 3){242.308}}
\put(2026,-1786){$x_8$}
\put(2626,-1336){$x_1$}
\put(3076,-886){$x_2$}
\put(3076,-361){$x_3$}
\put(2326,-211){$x_4$}
\put(1876,-436){$x_5$}
\put(1576,-1111){$x_6$}
\put(1576,-1636){$x_7$}
\put(2776,-2011){$y_1$}
\put(3676,-2011){$y_2$}
\put(1576,-2011){$z_1$}
\put(751,-2011){$z_2$}
\put(2326, 14){$A_4$}
\put(1201,-2311){$C_1$}
\put(3676,-1336){$a_2$}
\put(3601,-1636){$B_2$}
\put(3076,-1486){$a_1$}
\put(976,-1336){$a_6$}
\put(1276,-1786){$a_7$}
\put(2926,-1711){$B_1$}
\put(1801,-1861){\line( 1, 0){600}}
\put(2401,-1861){\line( 1, 1){900}}
\put(3301,-961){\line( 0, 1){300}}
\put(3301,-661){\line(-1, 1){600}}
\put(2701,-61){\line(-1, 0){600}}
\put(2101,-61){\line(-1,-1){600}}
\put(1501,-661){\line( 0,-1){900}}
\put(1501,-1561){\line( 1,-1){300}}
\end{picture}
\end{center}
\caption{}\label{Fig:ABC}
\end{figure}

Suppose that $n$ is odd and $n=2m+1$.  
By the sine rule, 
$$ \frac{a_j}{\sin j\theta_n} = \frac{y_j}{\sin \theta_n}
                              = \frac{a_{j-1}+x_j}{\sin(j+1)\theta_n} $$ 
for $0 < j \leq m-1$
where $\theta_n =2\pi/n$. Thus we have
$$ y_j\sin(j+1)\theta_n - y_{j-1}\sin(j-1)\theta_n = x_j\sin\theta_n. $$
The same argument for $C_j$ shows that
$$ z_j\sin(j+1)\theta_n - z_{j-1}\sin(j-1)\theta_n = x_{n-j}\sin\theta_n. $$
If we set
$$ w_0 := x_n + \sum_j y_j + \sum _j z_j, $$
the above equalities show that the coordinate
$(z_{m-1},\ldots,z_1,w_0,y_1,\ldots,y_{m-1})$ is obtained from
$(x_1,\ldots,x_n)$ by a linear isomorphism.

Since the areas of the triangles $A, B_j, C_j$ are as follows:
\begin{align*}
\Area A   & = w_0^2\frac14\tan\frac{\theta_n}2,\\
\Area B_j & = \frac12y_j(a_{j-1}+x_j)\sin j\theta_n
    = y_j^2 \, \frac{\sin j\theta_n\sin(j+1)\theta_n}{2\sin\theta_n}, \\
\Area C_j & = \frac12z_j(a_{n-j+1}+x_{n-j})\sin j\theta_n
    = z_j^2 \, \frac{\sin(-j)\theta_n\sin(-j+1)\theta_n}{2\sin\theta_n},
\end{align*}
we define the new coordinate by setting 
\begin{align*}
X_0 & = w_0\sqrt{\frac14\tan\frac{\theta_n}2}, \\
Y_j & = y_j\sqrt{\frac{\sin j\theta_n\sin(j+1)\theta_n}{2\sin\theta_n}}, \\
Z_j & = z_j\sqrt{\frac{\sin(-j)\theta_n\sin(-j+1)\theta_n}{2\sin\theta_n}}.
\end{align*}
$(Z_{m-1},\ldots,Z_1,X_0,Y_1,\ldots,Y_{m-1})$ is obtained from
$(z_{m-1},\ldots,z_1,w_0,y_1,\ldots,y_{m-1})$ by a 
linear isomorphism.
The area of the polygon is
$$ \Area = X_0^2 - \sum_j Y_j^2 - \sum_j Z_j^2. $$

If $n$ is even and $n = 2m + 2$, add $X_{m+1} :=
x_{m+1}\sqrt{(1/4)\tan\theta_n/2}$ and 
$(Z_{m-1},\ldots,Z_1,$ $X_0,Y_1,\ldots,Y_{m-1},X_{m+1})$ is the
coordinate.  
In this case, the area is 
$$ \Area = X_0^2 - X_{m+1}^2 - \sum_j Y_j^2 - \sum_j Z_j^2. $$

The function $\Area:{\cal E}_n^+ \rightarrow \real$ defines a
quadratic form of ${\cal E}_n$ which we also denote by $\Area$
and the above calculation shows that its signature is $(1, n-3)$.
\end{pf}

${\cal E}_n$  together with  $\Area$  becomes 
a Minkowski space.  
Let ${\cal P}_n$ be the intersection 
$\Area^{-1}(1)\cap \{X_0>0\}$.  
Then ${\cal P}_n$ is the hyperbolic space and $U_1$ is 
canonically homeomorphic to
$$ 
\Area^{-1}(1)\cap {\cal E}_n^+ 
= \Area^{-1}(1)\cap \bigcap_{i=1}^n \{x_i>0\} \subset {\cal P}_n.
$$
The region is bounded by ${\cal P}_n\cap \{x_i=0\}$  
for  $i = 1, 2, \cdots, n$.
Since $\{x_i=0\}$ is a hyperplane containing the origin in  
the Minkowski space  ${\cal E}_n$,  the
intersection with the hyperboloid ${\cal P}_n$ is the hyperbolic
hyperplane. It implies that:

\begin{LemNum}
There is a canonical homeomorphism between  $U_1$  and 
an interior of the $(n-3)$-dimensional
hyperbolic polyhedron obtained by 
taking closure of  $\Area^{-1}(1) \cap {\cal E}_n^+$.
\end{LemNum}

We denote this closed hyperbolic polyhedron 
by $\Delta_n$ and its faces
corresponding to $\{x_i=0\}\cap \Delta_n$ by $F_i$.
To see some geometric properties of  $\Delta_n$, 
we recall a lemma in hyperbolic geometry.
The Lorentz bilinear form  $q( \; , \; )$  on the Minkowski 
space  ${\cal E}_n$  is defined 
by a quadratic form  $\Area$  by setting  
$$
q(x,y) = 
\frac{1}{2}(\Area(x+y)-\Area(x)-\Area(y)).
$$  

\begin{LemNum}[\cite{Thurston},{\bf Proposition 2.4.5}]\label{Lem:Distance}
\begin{enumerate} 
\renewcommand{\labelenumi}{(\arabic{enumi})}%
	\item 
	Let $p_1$ and $p_2$ be in hyperbola ${\cal P}_n$ and  
	$d(p_1,p_2)$ the hyperbolic
	distance between $p_1$ and $p_2$. Then 
	$$
	\cosh d(p_1,p_2) = q(p_1,p_2).
	$$
	\item 
	Let $n_1$ and $n_2$ be normal vectors of the hyperplanes 
	$H_1$ and $H_2$  in ${\cal P}_n$.  
		\begin{enumerate} 
		\item 
	If  $H_1$  and  $H_2$  intersect in  ${\cal P}_n$,  
	then their dihedral angle  $\angle(H_1,H_2)$  satisfies 
	$$
	\cos\angle(H_1,H_2) = 
	\frac{q(n_1,n_2)}{\sqrt{q(n_1,n_1)} \, \sqrt{q(n_2,n_2)}}.
	$$
		\item 
	If  $H_1$  and  $H_2$  do not intersect in  ${\cal P}_n$,  
	then their shortest distance  $d(H_1,H_2)$  satisfies 
	$$
	\cosh d(H_1,H_2) = 
	\frac{q(n_1,n_2)}{\sqrt{q(n_1,n_1)} \, \sqrt{q(n_2,n_2)}}.
	$$
	\end{enumerate} 
\end{enumerate} 
\end{LemNum}

Then we have 

\begin{LemNum}\label{Lem:Intersection}
The faces of $\Delta_n$ intersect as follows.
\begin{enumerate}
\renewcommand{\labelenumi}{(\arabic{enumi})}%
\item $|i-j|\geq 2 \Rightarrow F_i \perp F_j$, 
\item If $n = 5$ or $6$, then $F_j\cap F_{j+1} = \emptyset$, 
\item If $n\geq 7$, 
$$ \cos(\omega_n) = \frac1{2\cos\frac{2\pi}n}, $$ 
where $\omega_n$ is the dihedral angle between $F_j$ and $F_{j+1}$.
\end{enumerate}
\end{LemNum}

\begin{pf}

(1) Suppose that $|i-j|\geq 2$.
We change the way we extend the edges of $P$ to form triangles $B_i$
and $C_i$ as we have done to get the new coordinate of ${\cal E}_n$.  
Extend two edges adjacent to the edge $x_i$ and $x_j$ respectively so
that we obtain an $(n - 2)$-gon $P^\prime$  
(see Figure \ref{Fig:Orthogonal}).  
\begin{figure}[ht]
\begin{center}
\setlength{\unitlength}{0.00083300in}%
\begingroup\makeatletter\ifx\SetFigFont\undefined%
\gdef\SetFigFont#1#2#3#4#5{%
  \reset@font\fontsize{#1}{#2pt}%
  \fontfamily{#3}\fontseries{#4}\fontshape{#5}%
  \selectfont}%
\fi\endgroup%
\begin{picture}(3924,1824)(289,-1873)
\put(2326,-286){$B_1$}
\thicklines
\put(1876,-961){\vector( 1, 0){300}}
\put(2551,-1861){\line( 1, 0){600}}
\put(3151,-1861){\line( 1, 1){600}}
\put(3751,-1261){\line( 0, 1){900}}
\put(3751,-361){\line(-1, 1){300}}
\put(3451,-61){\line(-1, 0){600}}
\put(2851,-61){\line(-1,-1){600}}
\put(2251,-661){\line( 0,-1){900}}
\put(2251,-1561){\line( 1,-1){300}}
\put(2851,-61){\line(-1, 0){600}}
\put(2251,-61){\line( 0,-1){600}}
\put(3751,-361){\line( 1,-1){450}}
\put(4201,-811){\line(-1,-1){450}}
\put(3826,-886){$B_2$}
\put(601,-511){$x_j$}
\put(1576,-961){$x_i$}
\put(601,-1861){\line( 1, 0){600}}
\put(1201,-1861){\line( 1, 1){600}}
\put(1801,-1261){\line( 0, 1){900}}
\put(1801,-361){\line(-1, 1){300}}
\put(1501,-61){\line(-1, 0){600}}
\put(901,-61){\line(-1,-1){600}}
\put(301,-661){\line( 0,-1){900}}
\put(301,-1561){\line( 1,-1){300}}
\end{picture}
\end{center}
\caption{}\label{Fig:Orthogonal}
\end{figure}

Call the triangles $B_1, B_2$.  For any linear isomorphism
on ${\cal E}_n$ whose first and second new coordinate are $Y_1 =
\sqrt{\Area(B_1)}$ and $Y_2 = \sqrt{\Area(B_2)}$, the
quadratic form Area is written as:
$$ 
\Area
= \text{(quadratic form without $Y_1$ and $Y_2$)} - Y_1^2 - Y_2^2.
$$
$\{x_i=0\}$ and $\{x_j=0\}$ in this coordinate is $\{Y_1 = 0\}$ and
$\{Y_2 = 0\}$ respectively and we can choose normal vectors  
by  $n_1 = (1, 0,\ldots, 0)$ and $n_2 = (0, 1, 0,\ldots, 0)$  
respectively.  
Therefore  $q(n_i, n_j) = 0$ and $F_i$ is orthogonal to $F_j$  
by Lemma~\ref{Lem:Distance} (2).

(2) Suppose that $n=5$.  
The intersection of $F_5$ and $F_1$ must consist of
(``generalized'' and ``degenerate'') equiangular pentagons as
in the figure.
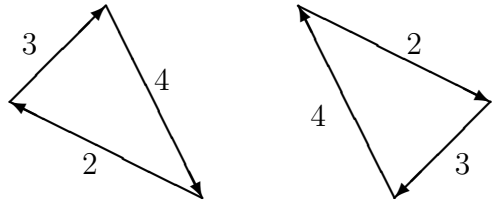
\begin{figure}[ht]
\begin{center}
\setlength{\unitlength}{0.00083300in}%
\begingroup\makeatletter\ifx\SetFigFont\undefined%
\gdef\SetFigFont#1#2#3#4#5{%
  \reset@font\fontsize{#1}{#2pt}%
  \fontfamily{#3}\fontseries{#4}\fontshape{#5}%
  \selectfont}%
\fi\endgroup%
\begin{picture}(3044,1244)(579,-683)
\put(3076,239){$2$}
\thicklines
\put(601,-61){\vector( 1, 1){600}}
\put(1201,539){\vector( 1,-2){600}}
\put(2401,539){\vector( 2,-1){1200}}
\put(3001,-661){\vector(-1, 2){600}}
\put(3601,-61){\vector(-1,-1){600}}
\put(676,239){$3$}
\put(1501, 14){$4$}
\put(1051,-511){$2$}
\put(2476,-211){$4$}
\put(3376,-511){$3$}
\put(1801,-661){\vector(-2, 1){1200}}
\end{picture}
\end{center}
\caption{Degenerated pentagons}
\end{figure}
There are two possibilities depending on the sign of $x_2$.  In both
cases the areas are negative 
so that they cannot intersect in ${\cal P}_5$.

The case $n=6$ is similar. If we set $x_i = x_{i+1} = 0$, then there
are four patterns of equiangular hexagons 
according to the sign of $x_{i-2}$ and $x_{i-1}$. 
The areas all are negative.

(3) By the symmetry of $\Delta_n$ we may assume that $j = 1$.
We denote $\{x_i = 0\}$ by $E_i$ for $i = 1, 2$.
$E_1$ is equal to $\{Y_1 = 0\}$ in the new coordinate.
The  $n$-gon with  $x_2 = 0$ is depicted as follows:
\begin{figure}[ht]
\begin{center}
\setlength{\unitlength}{0.00083300in}%
\begingroup\makeatletter\ifx\SetFigFont\undefined%
\gdef\SetFigFont#1#2#3#4#5{%
  \reset@font\fontsize{#1}{#2pt}%
  \fontfamily{#3}\fontseries{#4}\fontshape{#5}%
  \selectfont}%
\fi\endgroup%
\begin{picture}(5049,1776)(439,-2125)
\put(2176,-2086){$y_1$}
\thicklines
\put(3601,-1861){\line( 1, 4){375}}
\put(3901,-661){\line(-1,-1){1200}}
\put(3901,-661){\line(-2,-1){2400}}
\put(2401,-1261){$x_1$}
\put(3226,-1486){$y_1$}
\put(1951,-1786){$\theta_n$}
\put(3676,-1786){$3\theta_n$}
\put(2926,-1786){$2\theta_n$}
\put(3151,-2086){$y_2$}
\put(451,-1861){\line( 1, 0){5025}}
\end{picture}
\end{center}
\caption{}
\end{figure}

By the sine rule for the two triangles in the figure, we have
$$
\frac{y_1}{\sin 3\theta_n} = \frac{y_2}{\sin \theta_n}. 
$$
The ratio of the square of the areas of two triangles are
$$
\frac{Y_2^2}{Y_1^2} = \frac{y_2}{y_1} = \frac1{3-4\sin^2\theta_n}
$$
and constant. Thus
$$
E_2 = \{ Y_1 = Y_2\sqrt{3-4\sin^2\theta_n}\}. 
$$
Therefore normal vectors of these two hyperplanes can be 
chosen by 
\begin{eqnarray*}
	n_1 & = & (0,\ldots,0,1,0,\ldots,0),  \\
	n_2 & = & \left(0,\ldots,0,1,-\sqrt{3-4\sin^2\theta_n},\ldots,0\right).
\end{eqnarray*}
Then by Lemma~\ref{Lem:Distance} (1),  we have,  
$$ \cos \omega_n = \frac{q(n_1, n_2)}
        {\sqrt{q(n_1, n_1) \, q(n_2, n_2)} }
        = \frac1{2\cos\theta_n}. $$
\end{pf}

\begin{Remark}\label{Rem:angle}
{\em 
The formula in the Lemma tells us the angle.
\begin{center}\begin{tabular}{c|cccccc}
n                   & 7       & 8  & 9       & 10      & $\cdots$ &
$n\rightarrow\infty$ \\ \hline
$\omega_n$ (degree) & 36.6845 & 45 & 49.2542 & 51.8273 & $\cdots$ &
$\omega_n \nearrow 60$
\end{tabular}\end{center}
Therefore  $\omega_n<60$ for $n=7, 8, \ldots$.
}
\end{Remark}

\begin{Remark}
{\em 
$\Delta_5$ is bounded by 5 geodesics and $F_i$ and $F_{i\pm 2}$
are orthogonal, i.e. a hyperbolic right angle pentagon.
The same argument shows that $\Delta_6$ is a hyperbolic trigonal
bipyramid with angles between two intersecting faces $= \pi/2$
and three ideal vertices.
} 
\end{Remark}


\subsection{Identification of the boundary of  $X_n \approx X(n)$}

The number of hyperplanes $\{x_i=0\}$ with consecutive indices which
can meet on the boundary of $\Delta_n$ is at most $[\frac{n-1}2]-1$,
otherwise $[\frac{n-1}2]$  consecutive edges of 
the corresponding equiangular $n$-gon degenerate so
that the area cannot be positive. 
Such a degeneration corresponds to a collision 
of  $[\frac{n+1}2]$ consecutive points on the configuration 
which is ruled out in Yoshida's  $X_Y(n)$. 
Thus the boundary of $\Delta_n$ corresponds to
the set of degenerate configurations added in  $X_Y(n)$.  
On the other hand, 
our geometrization assigns not only 
a projective class of a configuration 
to each point in the interior of  $\Delta_n$, 
but a corresponding degenerate configuration to each point 
on the boundary of  $\Delta_n$.  
Gluing $(n-1)!/2$  copies of  
$\Delta_n$  together along the faces which 
represent the same degenerate configurations, 
and we obtain  $\overline{X_n}$  which is 
homeomorphic to Yoshida's compactification. 
$X_n \approx X(n)$  now lives in  
$\overline{X_n}$  as an open dense subset.  

We finish this section by discussing how singularities 
appear in  $\overline{X_n}$  by looking at 
gluing rule of polyhedral blocks.

The point lying on a face of codimension one in  $\Delta_n$   
corresponds to a configuration with a collision of two points.  
Hence we may label this face by a circular permutation of 
$n-2$ numbers and a group of two numbers up to reflection.  
For example,  the label  $(12) 3 \cdots n$  means that 
each point on that face represents a collision of the points 
marked by  $1$  and  $2$.   
Then the number of polyhedral blocks  
which share such a face is two according to how we approach 
to that degenerate configuration from nondegenerate ones. 
Hence the gluing does not yield any singularity along 
such a face.  
This proves that our geometrization 
is a hyperbolic cone-manifold by the definition of 
cone-manifolds.  

The point lying on a face of codimension two in  $\Delta_n$   
corresponds to a configuration with either 
a pair of collisions of two points  
or a collision of three point.  
We may label such a face by 
grouping marks involved in the collision together, such as 
$(12)3(45)6 \cdots n$,   
$(123)4 \cdots n$, etc. 

In the first case, 
the number of polyhedral blocks 
which share such a face is four according to how we approach 
to that degenerate configuration. 
On the other hand, 
the dihedral angle of two faces which share this codimension 
two face is  $\pi/2$  by Lemma~\ref{Lem:Intersection} (1).  
Hence again the gluing does not yield any singularity along 
such a face.  

These two observations show that  $\overline{X_n}$  is 
nonsingular when  $n = 5$.  
Actually, it is a hyperbolic surface which consists
of 12 hyperbolic right angle pentagons.  
Since each vertex belongs to four pentagons, 
the number of faces (pentagon), edges, vertices are 
$12, 30, 15$ respectively and Euler characteristic is  $-3$.  
It follows that it is a nonorientable surface homeomorphic to 
a connected sum of five copies of  
$\real \projective^2$.  

Let us discuss the case when  $n = 6$  before going into the other case. 
$X_6$ consists of the interior of 60 hyperbolic hexahedra.  
We are not allowed to have a collision of three successive points 
in this case.  
Hence the gluing does not yield any singularity 
along face of codimension at most two. 
Then consider a point on the face of codimension three.  
Since  $n -3 = 3$,  such a face is a vertex and 
corresponds to a triple of collisions of two points.  
The number of components of  $X_6$  
which share such a vertex is eight.  
On the other hand, the neighborhood of the vertex of $\Delta_6$ is
isometric to a neighborhood of the vertex of the intersection
between the non--negative orthant in the Poincare model 
of $\hyp^3$  in  $\real^3$.
Hence again the gluing does not yield any singularity.  
Moreover since horospherical cut of an ideal vertex 
in  $\Delta_6$  is always square, 
the gluing yields a complete end.  
Therefore  $\overline{X_6}$  is 
a complete hyperbolic 3-manifold.  

We can derive a few more information about geometry of  
$\overline{X_6}$.  
Since  $\Delta_6$  is scissors congruent to a quarter of 
the regular ideal octahedron, whose volume is  $3.66386 \cdots$,  
the volume of  $\overline{X_6}$  is  $54.957 \cdots$. 
$\overline{X_6}$  admits a natural action of the 
symmetry group of degree 6 by permuting labels of points.  
It turns out to be a full isometry group since 
the quotient is congruent to the smallest 
orbifold with appropriate date found by Adams in 
\cite{Adams}.  
We will see in the next section that  $\overline{X_6}$  has 
ten cusps.  

Suppose now that $n\geq 7$ and let us consider the neighborhood of a
degenerate configuration with a collision of three successive points.  
This always happens when  $n\geq 7$.  
There are six polyhedral blocks which share 
that configuration on the boundary according 
to the permutations of three numbers involved in the 
collision as in the figure below.  
\begin{figure}[ht]
\begin{center}
\setlength{\unitlength}{0.00083300in}%
\begingroup\makeatletter\ifx\SetFigFont\undefined%
\gdef\SetFigFont#1#2#3#4#5{%
  \reset@font\fontsize{#1}{#2pt}%
  \fontfamily{#3}\fontseries{#4}\fontshape{#5}%
  \selectfont}%
\fi\endgroup%
\begin{picture}(5424,2244)(589,-1723)
\put(3226,-1261){$3$}
\thicklines
\put(3601,-511){\circle*{100}}
\put(3901,-661){\circle*{100}}
\put(3301,-661){\circle*{100}}
\put(3301,-1411){\circle*{100}}
\put(3601,-1261){\circle*{100}}
\put(3901,-1411){\circle*{100}}
\put(5101,-1411){\circle*{100}}
\put(5401,-1261){\circle*{100}}
\put(5701,-1411){\circle*{100}}
\put(5701,-661){\circle*{100}}
\put(5401,-511){\circle*{100}}
\put(5101,-661){\circle*{100}}
\put(5101, 89){\circle*{100}}
\put(5401,239){\circle*{100}}
\put(5701, 89){\circle*{100}}
\put(3601,239){\circle*{100}}
\put(3901, 89){\circle*{100}}
\put(3301, 89){\circle*{100}}
\put(2101,-586){\vector( 1, 0){450}}
\put(601,-961){\line( 1, 1){600}}
\put(1201,-361){\line( 1,-1){600}}
\put(3001,-1711){\line( 1, 1){300}}
\put(3301,-1411){\line( 2, 1){300}}
\put(3601,-1261){\line( 2,-1){300}}
\put(3901,-1411){\line( 1,-1){300}}
\put(3001,-961){\line( 1, 1){300}}
\put(3301,-661){\line( 2, 1){300}}
\put(3601,-511){\line( 2,-1){300}}
\put(3901,-661){\line( 1,-1){300}}
\put(4801,-961){\line( 1, 1){300}}
\put(5101,-661){\line( 2, 1){300}}
\put(5401,-511){\line( 2,-1){300}}
\put(5701,-661){\line( 1,-1){300}}
\put(4801,-1711){\line( 1, 1){300}}
\put(5101,-1411){\line( 2, 1){300}}
\put(5401,-1261){\line( 2,-1){300}}
\put(5701,-1411){\line( 1,-1){300}}
\put(4801,-211){\line( 1, 1){300}}
\put(5101, 89){\line( 2, 1){300}}
\put(5401,239){\line( 2,-1){300}}
\put(5701, 89){\line( 1,-1){300}}
\put(3001,-211){\line( 1, 1){300}}
\put(3301, 89){\line( 2, 1){300}}
\put(3601,239){\line( 2,-1){300}}
\put(3901, 89){\line( 1,-1){300}}
\put(781,-151){$1 = 2 = 3$}
\put(3226,239){$1$}
\put(3526,389){$2$}
\put(3826,239){$3$}
\put(5026,239){$1$}
\put(5626,239){$2$}
\put(5326,389){$3$}
\put(3526,-361){$1$}
\put(3526,-1111){$1$}
\put(5626,-1261){$1$}
\put(5626,-511){$1$}
\put(3226,-511){$2$}
\put(5026,-511){$2$}
\put(3826,-1261){$2$}
\put(5326,-1111){$2$}
\put(3826,-511){$3$}
\put(5326,-361){$3$}
\put(5026,-1261){$3$}
\put(1201,-361){\circle*{100}}
\end{picture}
\end{center}
\end{figure}
But by Lemma~\ref{Lem:Intersection}, 
the dihedral angle of each piece is less than $2\pi/6$,
thus it gives rise to the singular points. 
Hence we have

\begin{ThNum}
$\overline{X_n}$  is a hyperbolic cone-manifold and 
homeomorphic to  $X_Y(n)$.
\begin{itemize}
\item When  $n = 5$ or $6$, it is nonsingular.
\item When  $n \geq 7$,  the singular set is nonempty. 
\end{itemize}
\end{ThNum}

\begin{Remark} 
{\em  
Every configuration appeared in  $\overline{X_n}$  has 
at least three labeled points which are disjointly 
placed on the circle.  
Normalize neighbor configurations of any particular 
degenerate one by sending such labeled points 
to  $\{0, 1, \infty \}$  by means of an projective automorphism.  
Then its neighborhood in  $\overline{X_n}$  is parameterized 
by the position of other  $n-3$  points.  
This observation shows that 
$\overline{X_n}$  is topologically a manifold.  
} 
\end{Remark}


\section{Deformation}


\subsection{Perturbation} 

In the previous section, the configuration space $X(n)$ is
identified with the space of marked equiangular $n$-gons up to
similarity.  
The identification is
given by the Schwarz--Christoffel mapping with all external angles
fixed to $2\pi/n$.  
If we perturb these external angles, the images of
Schwarz--Christoffel mapping change and one can expect that the
hyperbolic structure of the configuration space 
will deform accordingly.  

Let $\Theta_n$ be the set of $n$-tuples of real numbers $\theta =
(\theta_1, \ldots, \theta_n)$ satisfying the relations
$$ \sum_{j=1}^n\theta_j = 2\pi \quad \text{and} \quad
   0<\theta_i + \theta_j <\pi \quad (i,j\in \{1,\ldots, n\}).  $$
The indices should be understood modulo $n$ throughout the sequel.

Fix an element $\theta = (\theta_1, \ldots, \theta_n)$ of $\Theta_n$.
choose  $\alpha = (\alpha_1,\ldots, \alpha_n) 
	\in ({\real\projective}^1)^n - {\bold D}$,  
and we assign to  $\alpha$  
the unit disc in $\complex$ with  $n$  points specified on
the boundary. 
Then we map it conformally to an $n$-gon $P$ whose
vertices are the images of the specified points and the external angle
of the image of the $j$ th point  $\alpha_j$ is $\theta_j$. 
By the Schwarz--Christoffel formula,
$P$ is defined up to mark preserving similarity.  Let
$X_{n,\theta}$ be the space of mark preserving similarity classes
of Euclidean $n$-gons with external 
angles $\{\theta_1, \cdots, \theta_n\}$  compatible with 
markings. 
Then the same argument in
Lemma~\ref{Lem:Bijective} shows that

\begin{LemNum}\label{Lem:BijectiveWithTheta}
There is a canonical homeomorphism 
between  $X(n)$  and  $X_{n,\theta}$.
\end{LemNum}

As in the previous section, we get a hyperbolic polyhedral 
structure on each component
of $X_{n,\theta}$ --- we start from the parameterization
$(x_1,\ldots, x_n)$ of polygons by edge length and the function
$\Area$ turns out to be a quadratic form of signature $(1, n-3)$  again.  
Each component is identified with the intersection of hyperbola
$\Area^{-1}(1)$ and $n$ halfspaces.  
The components of  $X_{n, \theta}$  are no longer congruent 
each other. 
However the gluing rule still makes sense, and 
we obtain a hyperbolic cone-manifold  
$\overline{X_{n, \theta}}$  as well 
by identifying the boundary of  $X_{n, \theta} \approx X(n)$  
with corresponding weights. 
$\overline{X_{n, \theta}}$  contains  $X_{n, \theta}$  as 
an open dense subset. 
We will see what  $\overline{X_{n, \theta}}$  looks like 
when  $n = 5, 6$.


\subsection{Deformations of  $\overline{X_5}$}

Fix an element $\theta = (\theta_1, \ldots, \theta_5)$ of $\Theta_5$.
Let $p= \langle i_1 i_2 i_3 i_4 i_5 \rangle$  be a circular permutation of
$\{1,\ldots,5\}$ up to refection  
and $U_{p,\theta}$ a component of $X_{5,\theta}$
which consists of all pentagons whose 
marking correspond to  $p$. 
The external angle of the vertex marked
by  $i$ is $\theta_i$  by definition.  
Each element of $U_{p,\theta}$ is
parameterized by the side 
lengths  $(x_{i_1 i_2}, x_{i_2 i_3}, \ldots, x_{i_5 i_1})$  
where $x_{i_a i_b}$ is the length of the edge
between the vertices marked by  $i_a$ and $i_b$.

We modify the notation in the previous section and set
\begin{gather*}
{\cal E}_{p,\theta} := \left\{(x_{i_1 i_2},\ldots, x_{i_5 i_1})\mid
    x_{i_1 i_2}
   +x_{i_2 i_3}\exp(\sqrt{-1}\theta_{i_2})
    +\cdots+x_{i_5 i_1}\exp(\sqrt{-1}\sum_{j=2}^5\theta_{i_j})=0\right\},\\
X := \sqrt{\Area{A}}, \quad
Y := \sqrt{\Area{B}}, \quad
Z := \sqrt{\Area{C}}
\end{gather*}
where $A$, $B$ and $C$ are the triangles $A$, $B_1$ and $C_1$ in the
Figure~\ref{Fig:ABC} respectively.  By the same argument in the proof
of Lemma~\ref{Lem:Quadratic}, $(X, Y, Z)$ is a coordinate of 
${\cal E}_{p,\theta}$

Set ${\cal P}_5 = \Area^{-1}(1)\cap\{X>0\}$, then
$U_{p,\theta}$ is homeomorphic to  
${\cal P}_5\cap \bigcap_{a=1}^5 \{x_{i_a i_{a+1}}>0\}$ i.e. 
a hyperbolic pentagon.  
We denote by $\Delta_{p,\theta}$ 
this hyperbolic pentagon.  
We also simply denote by  $(i_1 i_2)i_3 i_4 i_5$  
its edge which corresponds to the degenerate
pentagons by the collision of the vertices $i_1$ and $i_2$.  
Similarly we use $(i_1 i_2)(i_3 i_4)i_5$, $(i_1 i_2 i_3)i_4 i_5$,  
etc, to represent faces of $\Delta_{p,\theta}$.

We next calculate the length of the edge $i_1 i_2 i_3(i_4 i_5)$.

\begin{LemNum}\label{Lem:EdgeLength}
Let $\theta=(\theta_1,\ldots,\theta_5)$ be an element of $\Theta_5$ and
$L(i_1 i_2 i_3(i_4 i_5);\theta)$ the length of the edge $i_1 i_2 i_3(i_4
i_5)$  of  $\Delta_{p,\theta}$.  
Then we have:
$$ \cosh L(i_1 i_2 i_3(i_4 i_5);\theta)
   = \sqrt{\frac{\sin\theta_{i_1}\sin\theta_{i_3}}
                {\sin(\theta_{i_1}+\theta_{i_2})
                 \sin(\theta_{i_2}+\theta_{i_3})}}. $$
\end{LemNum}

\begin{pf}
Suppose that the pentagon $\Delta_{p,\theta}$ is placed 
in the
hyperboloid of $XYZ$ space as explained above 
such that  $i_1 i_2 i_3(i_4i_5)$ and $i_1(i_2 i_3)i_4 i_5$  are 
in $\{Y=0\}$ and $\{Z=0\}$  respectively.
By Lemma~\ref{Lem:Intersection}, the end points of $i_1 i_2
i_3(i_4 i_5)$  are  $i_1(i_2 i_3)(i_4 i_5)$ and $(i_1 i_2)i_3(i_4 i_5)$.  
We denote their coordinates by $(x',y',z')$ and $(x'',y'',z'')$.

Observe that $i_1(i_2 i_3)(i_4 i_5)$ is in $\{Y=0\}$ and $\{Z=0\}$,
so $(x',y',z') = (1,0,0)$.

Next we calculate $x''$.  The pentagon $P$ with $x_{i_2 i_3}=0$ and
$x_{i_4 i_5}=0$ is depicted in the Figure~\ref{Fig:Intersect}.
\begin{figure}[ht]
\begin{center}
\setlength{\unitlength}{0.00083300in}%
\begingroup\makeatletter\ifx\SetFigFont\undefined%
\gdef\SetFigFont#1#2#3#4#5{%
  \reset@font\fontsize{#1}{#2pt}%
  \fontfamily{#3}\fontseries{#4}\fontshape{#5}%
  \selectfont}%
\fi\endgroup%
\begin{picture}(6324,1446)(-11,-1225)
\put(1726,-1186){$i_4$}
\thicklines
\multiput(4201, 89)(42.85714,-21.42857){8}{\makebox(6.6667,10.0000){\SetFigFont{10}{12}{\rmdefault}{\mddefault}{\updefault}.}}
\put(4501,-61){\line(-1,-1){900}}
\put(3601,-961){\line( 1, 0){2700}}
\put(6301,-961){\line(-2, 1){1800}}
\multiput(3301,-961)(42.85714,0.00000){8}{\makebox(6.6667,10.0000){\SetFigFont{10}{12}{\rmdefault}{\mddefault}{\updefault}.}}
\put(901,-961){\line(-1, 3){225}}
\put(1801,-961){\line( 1, 1){300}}
\multiput(601, 89)(42.85714,-21.42857){8}{\makebox(6.6667,10.0000){\SetFigFont{10}{12}{\rmdefault}{\mddefault}{\updefault}.}}
\multiput(901,-961)(15.00000,-45.00000){6}{\makebox(6.6667,10.0000){\SetFigFont{10}{12}{\rmdefault}{\mddefault}{\updefault}.}}
\put(1584,-701){\vector( 4,-1){300}}
\put(2401,-361){\vector( 1, 0){900}}
\put(443,-776){\vector( 4, 1){300}}
\put(901,-61){\line(-1,-1){900}}
\put(  1,-961){\line( 1, 0){2700}}
\put(2701,-961){\line(-2, 1){1800}}
\put(5101,-736){$P$}
\put(4426,-586){$e_2$}
\put(4351,-361){$\theta_{i_2}$}
\put(5326,-436){$e_1$}
\put(4501,-886){$\theta_{i_3}$}
\put(4651, 89){$A$}
\put(1201,-736){$P$}
\put(1201, 89){$A$}
\put(1051,-1186){$\theta_{i_3}$}
\put(4126,-136){$\theta_{i_1}$}
\put(526,-586){$\theta_{i_2}$}
\put(751,-1186){$i_3$}
\put(2101,-586){$i_5$}
\put(826, 14){$i_1$}
\put(751,-361){$i_2$}
\put(526,-136){$\theta_{i_1}$}
\put(2101,-916){$B$}
\put(3826,-511){$e_3$}
\put(2926,-886){$\theta_{i_2}+\theta_{i_3}$}
\put(4501,-61){\line( 1,-3){300}}
\end{picture}
\end{center}
\caption{Degeneration of $i_1 i_2$ and $i_4 i_5$}\label{Fig:Intersect}
\end{figure}
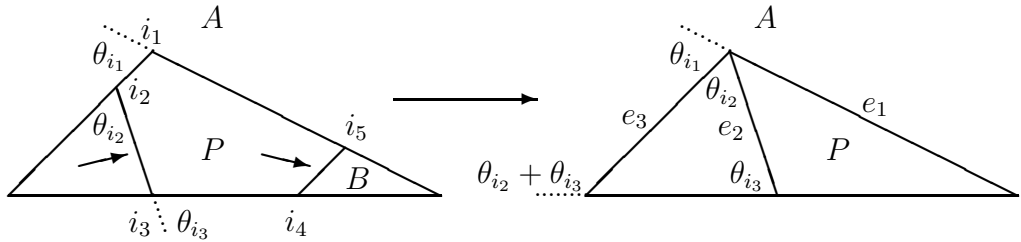
Calculate the areas of $A$ and $P$ by using the edge $e_1$ as the
common base edge:
$$\Area(P):\Area(A)
    = e_2\sin(\theta_{i_1}+\theta_{i_2}):e_3\sin\theta_{i_1}         $$
Since $\Area(P)=1$ at $(x'',y'',z'')$ and using the sine rule:
$$ \frac{e_2}{\sin(\theta_{i_2}+\theta_{i_3})}
    = \frac{e_3}{\sin\theta_{i_3}}                                   $$
we have
\begin{equation}\label{eq:x''}
x''
= \sqrt{\Area(A)}
= \sqrt{\frac{\sin\theta_{i_1}\sin\theta_{i_3}}
             {\sin(\theta_{i_1}+\theta_{i_2})\sin(\theta_{i_2}+\theta_{i_3})}}.
\end{equation}
By Lemma~\ref{Lem:Distance},
\begin{equation}\label{eq:coshl}
\cosh L(i_1 i_2 i_3(i_4 i_5);\theta)
   = q((x',y',z'),(x'',y'',z''))
   = q((1,0,0),(x'',y'',z'')) = x''.
\end{equation}
The identities (\ref{eq:x''}) and (\ref{eq:coshl}) prove the lemma.
\end{pf}

Now we investigate how the hyperbolic structure changes when $\theta$  
is perturbed. 
First, we comment on the topological
type of the new  $\overline{X_{5,\theta}}$.  
For each  $\theta$  near  
$\theta_0 = (2\pi/5, 2\pi/5, \cdots, 2\pi/5)$, 
the same argument in Lemma~\ref{Lem:Intersection} (1) shows that
$\Delta_{p,\theta}$  is a right angled pentagon, 
and hence  $\overline{X_{5,\theta}}$  is still a hyperbolic surface 
with the same topology as  
$\overline{X_5} \approx {\#}^5 \real\projective^2$.

The surface is nonorientable and does not support any  
complex structure at all.  
However we can still establish the analogue of 
Teichm\"uller theory.  
In fact, 
if we choose a maximal family of mutually disjoint 
nonparallel simple closed curves, then 
the set of hyperbolic structures is parameterized by 
their lengths and twisting amount for 2-sided ones.  
Hence the Teichm\"uller space  
${\cal T}({\#}^5 \real \projective^2)$  is homeomorphic 
to  $\real^9$.  
Following Teichm\"uller theory, 
we call this coordinate a Fenchel-Nielsen coordinate.  

To find a system of mutually disjoint, nonparallel simple closed
curves, let us enjoy some patch work.  
Here is part of our surface:
\begin{figure}[ht]
\begin{center}
\setlength{\unitlength}{0.00083300in}%
\begingroup\makeatletter\ifx\SetFigFont\undefined%
\gdef\SetFigFont#1#2#3#4#5{%
  \reset@font\fontsize{#1}{#2pt}%
  \fontfamily{#3}\fontseries{#4}\fontshape{#5}%
  \selectfont}%
\fi\endgroup%
\begin{picture}(4662,2724)(451,-4123)
\put(3076,-3736){\makebox(0,0)[lb]{\smash{\SetFigFont{12}{14.4}{\rmdefault}{\mddefault}{\updefault}213(45)}}}
\thicklines
\multiput(1351,-3211)(50.00000,0.00000){4}{\makebox(6.6667,10.0000){\SetFigFont{10}{12}{\rmdefault}{\mddefault}{\updefault}.}}
\multiput(1501,-3211)(40.00000,-20.00000){16}{\makebox(6.6667,10.0000){\SetFigFont{10}{12}{\rmdefault}{\mddefault}{\updefault}.}}
\multiput(2101,-3511)(0.00000,-50.00000){4}{\makebox(6.6667,10.0000){\SetFigFont{10}{12}{\rmdefault}{\mddefault}{\updefault}.}}
\multiput(3901,-1861)(0.00000,-50.00000){4}{\makebox(6.6667,10.0000){\SetFigFont{10}{12}{\rmdefault}{\mddefault}{\updefault}.}}
\multiput(3901,-2011)(40.00000,-20.00000){16}{\makebox(6.6667,10.0000){\SetFigFont{10}{12}{\rmdefault}{\mddefault}{\updefault}.}}
\multiput(4501,-2311)(50.00000,0.00000){4}{\makebox(6.6667,10.0000){\SetFigFont{10}{12}{\rmdefault}{\mddefault}{\updefault}.}}
\multiput(1351,-2761)(45.20548,0.00000){74}{\makebox(6.6667,10.0000){\SetFigFont{10}{12}{\rmdefault}{\mddefault}{\updefault}.}}
\multiput(2101,-1861)(0.00000,-50.00000){4}{\makebox(6.6667,10.0000){\SetFigFont{10}{12}{\rmdefault}{\mddefault}{\updefault}.}}
\multiput(2101,-2011)(-40.00000,-20.00000){16}{\makebox(6.6667,10.0000){\SetFigFont{10}{12}{\rmdefault}{\mddefault}{\updefault}.}}
\multiput(1501,-2311)(-50.00000,0.00000){4}{\makebox(6.6667,10.0000){\SetFigFont{10}{12}{\rmdefault}{\mddefault}{\updefault}.}}
\multiput(4651,-3211)(-50.00000,0.00000){4}{\makebox(6.6667,10.0000){\SetFigFont{10}{12}{\rmdefault}{\mddefault}{\updefault}.}}
\multiput(4501,-3211)(-40.00000,-20.00000){16}{\makebox(6.6667,10.0000){\SetFigFont{10}{12}{\rmdefault}{\mddefault}{\updefault}.}}
\multiput(3901,-3511)(0.00000,-50.00000){4}{\makebox(6.6667,10.0000){\SetFigFont{10}{12}{\rmdefault}{\mddefault}{\updefault}.}}
\put(901,-1711){\line( 1,-1){600}}
\put(1501,-2311){\line( 0,-1){900}}
\put(2101,-3511){\line( 1, 0){1800}}
\put(3901,-3511){\line( 1,-1){600}}
\put(1501,-1411){\line( 1,-1){600}}
\put(2101,-2011){\line( 1, 0){1800}}
\put(4501,-2311){\line( 0,-1){900}}
\put(4501,-3211){\line( 1,-1){600}}
\multiput(1501,-3211)(-30.00000,-30.00000){6}{\makebox(6.6667,10.0000){\SetFigFont{10}{12}{\rmdefault}{\mddefault}{\updefault}.}}
\multiput(2101,-3511)(-30.00000,-30.00000){6}{\makebox(6.6667,10.0000){\SetFigFont{10}{12}{\rmdefault}{\mddefault}{\updefault}.}}
\multiput(4051,-1861)(-30.00000,-30.00000){6}{\makebox(6.6667,10.0000){\SetFigFont{10}{12}{\rmdefault}{\mddefault}{\updefault}.}}
\multiput(4651,-2161)(-30.00000,-30.00000){6}{\makebox(6.6667,10.0000){\SetFigFont{10}{12}{\rmdefault}{\mddefault}{\updefault}.}}
\put(2251,-1936){\makebox(0,0)[lb]{\smash{\SetFigFont{12}{14.4}{\rmdefault}{\mddefault}{\updefault}124(35)}}}
\put(3076,-1936){\makebox(0,0)[lb]{\smash{\SetFigFont{12}{14.4}{\rmdefault}{\mddefault}{\updefault}214(35)}}}
\put(2176,-2461){\makebox(0,0)[lb]{\smash{\SetFigFont{12}{14.4}{\rmdefault}{\mddefault}{\updefault}12435}}}
\put(3451,-2461){\makebox(0,0)[lb]{\smash{\SetFigFont{12}{14.4}{\rmdefault}{\mddefault}{\updefault}21435}}}
\put(4576,-2611){\makebox(0,0)[lb]{\smash{\SetFigFont{12}{14.4}{\rmdefault}{\mddefault}{\updefault}143(52)}}}
\put(826,-2611){\makebox(0,0)[lb]{\smash{\SetFigFont{12}{14.4}{\rmdefault}{\mddefault}{\updefault}243(51)}}}
\put(826,-3061){\makebox(0,0)[lb]{\smash{\SetFigFont{12}{14.4}{\rmdefault}{\mddefault}{\updefault}234(51)}}}
\put(2176,-3211){\makebox(0,0)[lb]{\smash{\SetFigFont{12}{14.4}{\rmdefault}{\mddefault}{\updefault}12345}}}
\put(3451,-3211){\makebox(0,0)[lb]{\smash{\SetFigFont{12}{14.4}{\rmdefault}{\mddefault}{\updefault}21345}}}
\put(4201,-3661){\makebox(0,0)[lb]{\smash{\SetFigFont{12}{14.4}{\rmdefault}{\mddefault}{\updefault}23145}}}
\put(4951,-3586){\makebox(0,0)[lb]{\smash{\SetFigFont{12}{14.4}{\rmdefault}{\mddefault}{\updefault}314(52)}}}
\put(4576,-3061){\makebox(0,0)[lb]{\smash{\SetFigFont{12}{14.4}{\rmdefault}{\mddefault}{\updefault}134(52)}}}
\put(2251,-3736){\makebox(0,0)[lb]{\smash{\SetFigFont{12}{14.4}{\rmdefault}{\mddefault}{\updefault}123(45)}}}
\put(451,-2086){\makebox(0,0)[lb]{\smash{\SetFigFont{12}{14.4}{\rmdefault}{\mddefault}{\updefault}423(51)}}}
\put(1276,-1936){\makebox(0,0)[lb]{\smash{\SetFigFont{12}{14.4}{\rmdefault}{\mddefault}{\updefault}14235}}}
\put(1726,-1561){\makebox(0,0)[lb]{\smash{\SetFigFont{12}{14.4}{\rmdefault}{\mddefault}{\updefault}142(35)}}}
\put(3601,-4036){\makebox(0,0)[lb]{\smash{\SetFigFont{12}{14.4}{\rmdefault}{\mddefault}{\updefault}231(45)}}}
\multiput(3001,-3661)(0.00000,45.00000){41}{\makebox(6.6667,10.0000){\SetFigFont{10}{12}{\rmdefault}{\mddefault}{\updefault}.}}
\end{picture}
\end{center}
\caption{}\label{Fig:Patch}
\end{figure}

The end points of each solid
lines in the figure are identified in $\overline{X_{5, \theta}}$  to 
form closed geodesics and they are mutually disjoint. 
Reading the labels in the figure, we have simple closed curves
$(i_4 i_5)$ in $\overline{X_{5,\theta}}$ which consists of 
three edges
$i_1 i_2 i_3(i_4 i_5)$, $i_2 i_1 i_3(i_4 i_5)$  and  $i_2 i_3 i_1(i_4 i_5)$. 

We choose $(j5)$, $j=1, \ldots,4$ as a
system of four mutually disjoint, nonparallel simple closed curves, 
and let $L(ij;\theta)$ denote the length of $(ij)$ where $i,j\in
\{1,\ldots,5\}$ and $\theta\in\Theta_5$.  In the next lemma, we
calculate  $L(ij;\theta)$.  
For $\theta=(\theta_1, \theta_2, \theta_3, \theta_4, \theta_5)\in\Theta$, 
define a function
$N(i_1i_2i_3;\theta)$ by
\begin{align*}
&N(i_1i_2i_3;\theta) \\
&=\frac{\sin\theta_{i_1}\sin\theta_{i_2}\sin\theta_{i_3}
         -\sin(\theta_{i_1}\!+\!\theta_{i_2}\!+\!\theta_{i_3})
          (\sin\theta_{i_1}\sin\theta_{i_2}
	  +\sin\theta_{i_2}\sin\theta_{i_3}
          +\sin\theta_{i_3}\sin\theta_{i_1})}
        {\sin(\theta_{i_1}\!+\theta_{i_2})
         \sin(\theta_{i_2}\!+\theta_{i_3})
         \sin(\theta_{i_3}\!+\theta_{i_1})}.
\end{align*}
Then we have:

\begin{LemNum}\label{Lem:CurveLength}
For $\{i_1,i_2,i_3,i_4,i_5\}=\{1,\ldots,5\}$ and $\theta =
(\theta_i, \theta_2, \theta_3, \theta_4, \theta_5)
\in \Theta_5$:
$$ \cosh(L(i_4i_5;\theta)) = N(i_1i_2i_3;\theta). $$
\end{LemNum}

To prove this lemma, we prepare a lemma for trigonometric computation.

\begin{LemNum}
For any $\alpha, \beta, \gamma\in \real$, we have
\begin{equation}
\sin\alpha\sin\gamma - \sin(\alpha+\beta)\sin(\beta+\gamma)
	= -\sin\beta\sin(\alpha+\beta+\gamma).           \label{eq:sin}
\end{equation}
\end{LemNum}

\begin{pf}
The addition rule says the identity,  
$(\cos(a-b)-\cos(a+b))/2 = \sin a \sin b$. 
Apply
this relation for each term in the equation~(\ref{eq:sin}).
\end{pf}

\begin{pf*}{\it Proof of Lemma~\ref{Lem:CurveLength}.}
Put $\Delta_{p,\theta}$ in the hyperbola 
of $XYZ$ space again (see
Figure~\ref{Fig:Develope}(a)).
\begin{figure}[ht]
\begin{center}
\setlength{\unitlength}{0.00083300in}%
\begingroup\makeatletter\ifx\SetFigFont\undefined%
\gdef\SetFigFont#1#2#3#4#5{%
  \reset@font\fontsize{#1}{#2pt}%
  \fontfamily{#3}\fontseries{#4}\fontshape{#5}%
  \selectfont}%
\fi\endgroup%
\begin{picture}(5874,2373)(139,-1972)
\put(2776,-436){$\phi$}
\thicklines
\put(4201,-61){\line( 0,-1){450}}
\put(4201,-811){\line( 0,-1){150}}
\put(5101,-811){\line( 0,-1){150}}
\put(601,-61){\line( 1, 0){1800}}
\put(2401,-61){\line( 0,-1){900}}
\put(2401,-961){\line(-3,-2){450}}
\put(1951,-1261){\line(-3, 2){450}}
\put(1501,-961){\line(-3,-2){450}}
\put(1051,-1261){\line(-3, 2){450}}
\put(601,-961){\line( 0, 1){900}}
\put(1501,-61){\line( 0,-1){450}}
\put(1501,-811){\line( 0,-1){150}}
\put(1501,-1411){\vector( 0,-1){0}}
\put(5101,-61){\line( 0,-1){450}}
\put(151,-61){\vector(-1, 0){0}}
\put(2851,-61){\vector(-1, 0){0}}
\put(5101,-1411){\vector( 0,-1){0}}
\put(2551,-511){\vector( 1, 0){600}}
\put(4276,-361){$e_2$}
\put(3976,-736){(12)345}
\put(826,-511){12345}
\put(1201,-736){(12)345}
\put(976,-211){$e_1$}
\put(1576,-361){$e_2$}
\put(1876,-211){$e_3$}
\put(826, 14){123(45)}
\put(1651, 14){213(45)}
\put(3526, 14){123(45)}
\put(3526,-511){12345}
\put(3676,-211){$e_1$}
\put(4576,-211){$e_3$}
\put(4426, 14){213(45)}
\put(5326, 14){231(45)}
\put(5476,-211){$e_4$}
\put(5401,-511){23145}
\put(4801,-736){2(13)45}
\put(2851, 14){$Y$}
\put(151, 14){$Y$}
\put(526, 14){$p_1$}
\put(6001, 14){$p_4$}
\put(1801,-511){21345}
\put(4501,-511){21345}
\put(1426,-1636){$Z$}
\put(5026,-1636){$Z$}
\put(2401, 14){$p_3$}
\put(3076, 14){$\phi(p_1)$}
\put(1426,-1936){(a)}
\put(4576,-1936){(b)}
\put(3301,-61){\line( 0,-1){900}}
\put(3301,-961){\line( 3,-2){450}}
\put(3751,-1261){\line( 3, 2){450}}
\put(4201,-961){\line( 3,-2){450}}
\put(4651,-1261){\line( 3, 2){450}}
\put(5101,-961){\line( 3,-2){450}}
\put(5551,-1261){\line( 3, 2){450}}
\put(6001,-961){\line( 0, 1){900}}
\put(6001,-61){\line(-1, 0){2700}}
\multiput(1501,389)(0.00000,-45.00000){11}{\makebox(6.6667,10.0000){\SetFigFont{10}{12}{\rmdefault}{\mddefault}{\updefault}.}}
\multiput(1501,-961)(0.00000,-45.00000){10}{\makebox(6.6667,10.0000){\SetFigFont{10}{12}{\rmdefault}{\mddefault}{\updefault}.}}
\multiput(601,-61)(-45.00000,0.00000){10}{\makebox(6.6667,10.0000){\SetFigFont{10}{12}{\rmdefault}{\mddefault}{\updefault}.}}
\multiput(3301,-61)(-45.00000,0.00000){10}{\makebox(6.6667,10.0000){\SetFigFont{10}{12}{\rmdefault}{\mddefault}{\updefault}.}}
\multiput(5101,-961)(0.00000,-45.00000){10}{\makebox(6.6667,10.0000){\SetFigFont{10}{12}{\rmdefault}{\mddefault}{\updefault}.}}
\multiput(5101,389)(0.00000,-45.00000){11}{\makebox(6.6667,10.0000){\SetFigFont{10}{12}{\rmdefault}{\mddefault}{\updefault}.}}
\end{picture}
\end{center}
\caption{}\label{Fig:Develope}
\end{figure}
Suppose that $e_1 := i_1 i_2 i_3(i_4 i_5)$
corresponds to $\{Y=0\}$ and $e_2 := (i_1 i_2) i_3 i_4 i_5$
corresponds to $\{Z=0\}$.  Developing across $e_2$, we meet $i_2 i_1
i_3 i_4 i_5$ and the edge next to $e_1$ is $e_3 := i_2 i_1 i_3(i_4
i_5)$. Let $p_a = (x_a, y_a, z_a)$ be the end points of $e_a$ $(a =
1,3)$ which is not $(1,0,0)$. For simplicity, we denote
$\sin\theta_{i_j}$, $\sin(\theta_{i_j}+\theta_{i_k})$,
$\sin(\theta_{i_j}+\theta_{i_k}+\theta_{i_l})$ by $s_j$, $s_{jk}$,
$s_{jkl}$ respectively.  Then by Lemma~\ref{Lem:EdgeLength}
\begin{alignat*}{3}
x_1 & =  \sqrt{\frac{s_1 s_3} {s_{12}s_{23}}},&\quad
y_1 & =  \sqrt{{x_1}^2 -1}
      =  \sqrt{\frac{s_1 s_3 - s_{12}s_{23}}{s_{12}s_{23}}}
      =  \sqrt{\frac{-s_2 s_{123}}{s_{12}s_{23}}},&\quad
z_1 & =  0, \\
x_3 & =  \sqrt{\frac{s_2 s_3} {s_{21}s_{13}}},&\quad
y_3 & = -\sqrt{{x_3}^2 -1}
      = -\sqrt{\frac{s_2 s_3 - s_{21}s_{13}}{s_{21}s_{13}}}
      = -\sqrt{\frac{-s_1 s_{123}}{s_{21}s_{13}}},&\quad
z_3 & = 0. 
\end{alignat*}
We used equation (\ref{eq:sin}) for $y_1$ and $y_3$.

Now apply the hyperbolic isometry $\phi$ on the hyperboloid which
fixes the geodesic through $e_1$ (and $e_3$) setwise and sends $p_3$
to $(1,0,0)$ (Figure~\ref{Fig:Develope}(b)).  Denote by $(x_1', y_1',
z_1')$ the coordinate of $\phi(p_1)$. Since $\phi$ preserves $q$
\begin{align*}
x'_1 & = x_1'\cdot 1 - y_1'\cdot 0 - z_1'\cdot 0 
       = q(\phi(p_1), \phi(p_3)) = q(p_1, p_3) 
       = x_1 x_3 - y_1 y_3 \\
     & = \sqrt{\frac{s_1 s_3 s_2 s_3}{s_{12}s_{23}s_{21}s_{13}}}
         + \sqrt{\frac{s_2 s_{123} s_1 s_{123}}{s_{12}s_{23}s_{21}s_{13}}}, \\
y'_1 & = \sqrt{{x'_1}^2 -1}
       = \sqrt{\frac{s_1 s_3 s_2 s_3 + s_2 s_{123}s_1 s_{123}
                     - 2  s_1 s_2 s_3 s_{123} - s_{12}s_{23}s_{21}s_{13}}
                    {s_{12}s_{23}s_{21}s_{13}}} \\
\intertext{In the above equation, we used $\sqrt{s_i^2}=s_i,
\sqrt{s_{ij}^2}=s_{ij}$ and $\sqrt{s_{ijk}^2}=-s_{ijk}$ because
$\theta_{i_j}+\theta_{i_{j+1}}<\pi$.
By (\ref{eq:sin}), $s_{12}s_{23}s_{21}s_{13}=(s_1 s_3+s_2
s_{123})(s_2 s_3+s_1 s_{123})$, hence}
     & = \sqrt{\frac{s_1 s_3 s_2 s_3 + s_2 s_{123}s_1 s_{123}
                     - 2  s_1 s_2 s_3 s_{123}
                     - (s_1 s_3+s_2 s_{123})(s_2 s_3+s_1 s_{123})}
                    {s_{12}s_{23}s_{21}s_{13}}} \\
     & = \sqrt{\frac{-(s_1 + s_2)^2 s_3 s_{123}}
                    {s_{12}s_{23}s_{21}s_{13}}}, \\
z'_1 & = 0. 
\end{align*}
Develop again across $i_2(i_1 i_3)i_4 i_5$ to meet the component
$i_2 i_3 i_1 i_4 i_5$ and the edge next to $e_3$ is  $e_4 = i_2 i_3
i_1(i_4 i_5)$. Let $p_4 = (x_4, y_4, z_4)$ be the end points of $e_4$
which is not $(1,0,0)$. Then again by lemma \ref{Lem:EdgeLength}
\begin{alignat*}{3}
x_4 & =  \sqrt{\frac{s_2 s_1} {s_{23}s_{31}}}, &\quad
y_4 & = -\sqrt{{x_4}^2 -1}
      = -\sqrt{\frac{s_2 s_1 - s_{23}s_{31}}{s_{23}s_{31}}}
      = -\sqrt{\frac{-s_3 s_{123}}{s_{23}s_{31}}}, &\quad
z_4 & =  0. 
\end{alignat*}
Hence, the hyperbolic cosine of the length of $e_1 e_3 e_4$ is 
\begin{align*}
\cosh & d(\phi(p_1),p_4)  = x'_1 x_4 - y'_1 y_4 - z'_1 z_4\\
 & = (\sqrt{\frac{s_1 s_3 s_2 s_3}{s_{12}s_{23}s_{21}s_{13}}}
         + \sqrt{\frac{s_2 s_{123} s_1s_{123}}{s_{12}s_{23}s_{21}s_{13}}} )
     \sqrt{\frac{s_2 s_1} {s_{23}s_{31}}}
     + \sqrt{\frac{-(s_1 + s_2)^2 s_3 s_{123}}{s_{12}s_{23}s_{21}s_{13}}}
       \sqrt{\frac{-s_3 s_{123}}{s_{23}s_{31}}} \\
 & = \frac{s_1 s_2 s_3 - s_{123}(s_1 s_2 + s_2 s_3 + s_3 s_1)}
          {s_{12}s_{23}s_{31}}.   
\end{align*}
\end{pf*}


\begin{ThNum}\label{Thm:DeformationOfX(5)}
The map from $\Theta_5$ to  ${\cal T}({\#}^5 \real \projective^2)$  
defined by the assignment
$: \theta \mapsto \overline{X_{n,\theta}}$ is a local embedding 
at  $\theta_0=(2\pi/5,\ldots,2\pi/5)$.
\end{ThNum}

\begin{pf}
Choose four closed geodesics $(15), (25), (35), (45)$.  
By Figure~\ref{Fig:Patch}, they are mutually disjoint and nonparallel, 
so that their lengths will be 
a part of a Fenchel-Nielsen coordinate.  
Taking these components of the map from  $\Theta$  to  
${\cal T}({\#}^5 \real \projective^2)$, we get  
\begin{alignat*}{4}
&\Phi: &                   &\Theta_5          && \rightarrow && \real^4 \\
&      & \theta=(\theta_1,& \ldots, \theta_5) && \mapsto     &
(L(15;\theta),L(25;\theta),&L(35;\theta),L(45;\theta)).
\end{alignat*}
It is enough to show that 
the Jacobian of $\Phi$ does not vanish at $\theta_0$.
As a basis of the tangent space of $\Theta$ at $\theta_0$, we take
four paths
$p_j(t) = (\theta_{j1}(t), \cdots, \theta_{j5}(t))$ 
$j=1,\ldots,4$ passing through the barycenter $\theta_0$ defined by
$\theta_{jj}(t)=\frac{2\pi}5 +t$,
$\theta_{jj+1}(t)=\frac{2\pi}5 -t$,
$\theta_{jk}(t)=\frac{2\pi}5$ for $k \ne j, j+1$.
Then,
\begin{align*}
\Phi(p_1(t)) & =(a(t), b(t), c(t), c(t)), \\
\Phi(p_2(t)) & =(c(t), a(t), b(t), c(t)), \\
\Phi(p_3(t)) & =(c(t), c(t), a(t), b(t)), \\
\Phi(p_4(t)) & =(b(t), b(t), b(t), d(t)), 
\end{align*}
where
\begin{alignat*}{2}
\cosh a(t) &= N(234;p_1(t)),&\quad
\cosh b(t) &= N(134;p_1(t)),\\
\cosh c(t) &= N(124;p_1(t)),&\quad
\cosh d(t) &= N(123;p_4(t)).
\end{alignat*}
Then since 
$a'(0)=-b'(0)$ is non zero and 
$c'(0)=d'(0)=0$,
it is easy to see that the Jacobian of $\Phi$ at $\theta_{0}$ 
does not vanish.
\end{pf}


\subsection{Deformations of  $\overline{X_6}$}

As in the previous case, we use the notation $\Theta_6$, 
$p= \langle i_1 i_2 i_3 i_4 i_5 i_6 \rangle$,  
$U_{p,\theta}$ and ${\cal P}_6$, etc.

Also set
\begin{gather*}
{\cal E}_{p,\theta} := \left\{(x_{i_1 i_2},\ldots, x_{i_6 i_1})\mid
    x_{i_1 i_2}
    +x_{i_2 i_3}\exp(\sqrt{-1}\theta_{i_2})
    +\cdots+x_{i_6 i_1}\exp(\sqrt{-1}\sum_{j=2}^6\theta_{i_j})=0\right\},\\
X := \sqrt{\Area{A}}, \quad
Y := \sqrt{\Area{B}}, \quad
Z := \sqrt{\Area{C}}, \quad
W := \sqrt{\Area{D}}
\end{gather*}
where $A$,$B$,$C$,$D$ are  $A$,$B_1$,$C_1$,$A_3$ in
Figure~\ref{Fig:ABC}. Then $(X, Y, Z, W)$ is a coordinate of ${\cal
E}_{p,\theta}$.  $U_{p,\theta}$ is homeomorphic to the region
$\Delta_{p,\theta}={\cal P}_6\cap \bigcap_{a=1}^6 
\{x_{i_a i_{a+1}} \geq 0\}$  which is a hyperbolic hexahedron 
of finite volume.

\begin{figure}[ht]
\begin{center}
\setlength{\unitlength}{0.00083300in}%
\begingroup\makeatletter\ifx\SetFigFont\undefined%
\gdef\SetFigFont#1#2#3#4#5{%
  \reset@font\fontsize{#1}{#2pt}%
  \fontfamily{#3}\fontseries{#4}\fontshape{#5}%
  \selectfont}%
\fi\endgroup%
\begin{picture}(5274,2331)(589,-1948)
\put(5476,-436){\makebox(0,0)[lb]{\smash{\SetFigFont{12}{14.4}{\rmdefault}{\mddefault}{\updefault}(45)}}}
\thicklines
\put(1501,239){\line(-1,-4){300}}
\put(1201,-961){\line( 1,-2){300}}
\put(1501,-1561){\line( 1, 1){900}}
\put(601,-661){\line( 1,-1){900}}
\multiput(601,-661)(45.00000,0.00000){41}{\makebox(6.6667,10.0000){\SetFigFont{10}{12}{\rmdefault}{\mddefault}{\updefault}.}}
\put(826,-1336){\vector( 1, 1){300}}
\put(1501,-136){\vector(-1,-1){0}}
\multiput(1801,164)(-33.33333,-33.33333){9}{\makebox(6.6667,10.0000){\SetFigFont{10}{12}{\rmdefault}{\mddefault}{\updefault}.}}
\put(1801,-1036){\vector(-1, 1){0}}
\multiput(2101,-1336)(-33.33333,33.33333){9}{\makebox(6.6667,10.0000){\SetFigFont{10}{12}{\rmdefault}{\mddefault}{\updefault}.}}
\put(1351,-961){\vector( 1, 0){1650}}
\put(3601,-661){\line( 0,-1){600}}
\put(3301,-961){\line( 1, 0){600}}
\put(4201,-661){\vector( 4, 1){600}}
\put(4201,-1261){\vector( 4,-1){600}}
\put(5026,-361){\line( 1, 0){300}}
\put(5326,-361){\line( 0,-1){300}}
\put(5326,-361){\line( 1, 1){225}}
\put(5326,-1111){\line( 0,-1){300}}
\put(5326,-1411){\line(-1, 0){300}}
\put(5326,-1411){\line( 1,-1){225}}
\put(5551,164){\line( 0,-1){300}}
\put(5551,-136){\line( 1, 0){300}}
\put(5851,-1636){\line(-1, 0){300}}
\put(5551,-1636){\line( 0,-1){300}}
\put(1576,-436){\makebox(0,0)[lb]{\smash{\SetFigFont{12}{14.4}{\rmdefault}{\mddefault}{\updefault}(12)}}}
\put(1801,239){\makebox(0,0)[lb]{\smash{\SetFigFont{12}{14.4}{\rmdefault}{\mddefault}{\updefault}(34)}}}
\put(676,-1561){\makebox(0,0)[lb]{\smash{\SetFigFont{12}{14.4}{\rmdefault}{\mddefault}{\updefault}(23)}}}
\put(2026,-1561){\makebox(0,0)[lb]{\smash{\SetFigFont{12}{14.4}{\rmdefault}{\mddefault}{\updefault}(61)}}}
\put(3226,-1186){\makebox(0,0)[lb]{\smash{\SetFigFont{12}{14.4}{\rmdefault}{\mddefault}{\updefault}(23)}}}
\put(3676,-1186){\makebox(0,0)[lb]{\smash{\SetFigFont{12}{14.4}{\rmdefault}{\mddefault}{\updefault}(45)}}}
\put(3226,-886){\makebox(0,0)[lb]{\smash{\SetFigFont{12}{14.4}{\rmdefault}{\mddefault}{\updefault}(56)}}}
\put(3676,-886){\makebox(0,0)[lb]{\smash{\SetFigFont{12}{14.4}{\rmdefault}{\mddefault}{\updefault}(12)}}}
\put(4951,-1336){\makebox(0,0)[lb]{\smash{\SetFigFont{12}{14.4}{\rmdefault}{\mddefault}{\updefault}(56)}}}
\put(4951,-586){\makebox(0,0)[lb]{\smash{\SetFigFont{12}{14.4}{\rmdefault}{\mddefault}{\updefault}(23)}}}
\put(901,-511){\makebox(0,0)[lb]{\smash{\SetFigFont{12}{14.4}{\rmdefault}{\mddefault}{\updefault}(56)}}}
\put(1426,-1186){\makebox(0,0)[lb]{\smash{\SetFigFont{12}{14.4}{\rmdefault}{\mddefault}{\updefault}(45)}}}
\put(3451,-436){$\theta_4+\theta_5+\theta_6<\pi$}
\put(3451,-1636){$\theta_1+\theta_2+\theta_3<\pi$}
\put(5626,-61){\makebox(0,0)[lb]{\smash{\SetFigFont{12}{14.4}{\rmdefault}{\mddefault}{\updefault}(12)}}}
\put(5626,-1861){\makebox(0,0)[lb]{\smash{\SetFigFont{12}{14.4}{\rmdefault}{\mddefault}{\updefault}(45)}}}
\put(5101,-211){\makebox(0,0)[lb]{\smash{\SetFigFont{12}{14.4}{\rmdefault}{\mddefault}{\updefault}(56)}}}
\put(5476,-1411){\makebox(0,0)[lb]{\smash{\SetFigFont{12}{14.4}{\rmdefault}{\mddefault}{\updefault}(12)}}}
\put(5101,-1711){\makebox(0,0)[lb]{\smash{\SetFigFont{12}{14.4}{\rmdefault}{\mddefault}{\updefault}(23)}}}
\put(601,-661){\line( 2,-1){600}}
\put(1201,-961){\line( 4, 1){1200}}
\put(2401,-661){\line(-1, 1){900}}
\put(1501,239){\line(-1,-1){900}}
\end{picture}
\end{center}
\caption{Cusp of $\overline{X_6}$. }\label{Fig:CuspDeform}
\end{figure}
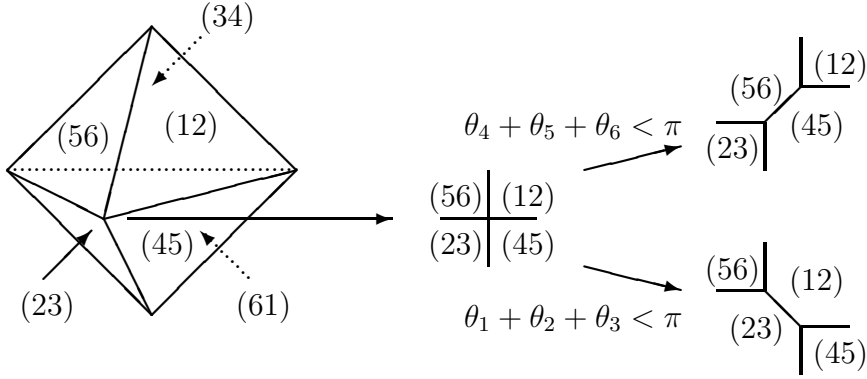

Let us describe how $\Delta_{p,\theta}$  deforms when we perturb
$\theta$ (see Figure~\ref{Fig:CuspDeform}).  We denote the
face of $\Delta_{p,\theta}$ which corresponds to the set of degenerate
hexagons by collisions of the points $i_j$ and $i_{j+1}$ by
$(i_ji_{j+1})i_{j+2} i_{j+3} i_{j+4} i_{j+5}$ or $(i_ji_{j+1})$ if
there is no confusion.  

When the weights are equal, namely  $\theta_0 = (2\pi/6, \cdots, 2\pi/6)$,
$\Delta_{p,\theta_0}$  has three ideal vertices.
Observe that the four faces containing an ideal vertex has labels of
type $(i_{k}i_{k+1})$, $(i_{k+1}i_{k+2})$, $(i_{k+3}i_{k+4})$,
$(i_{k+4}i_{k+5})$ for some $k\in\{0,\ldots,5\}$.  
We denote this
vertex by $(i_{k}i_{k+1}i_{k+2}) (i_{k+3}i_{k+4}i_{k+5})$.  

If we perturb the angles so that  $\theta_{i_k} + \theta_{i_{k+1}} +
\theta_{i_{k+2}} < \pi$, then three vertices $i_{k}$, $i_{k+1}$,
$i_{k+2}$ of the hexagon can collide, 
and  $(i_{k}i_{k+1})$ and $(i_{k+1}i_{k+2})$ intersects in ${\cal P}_6$.  
If
$\theta_{i_{k+3}} + \theta_{i_{k+4}} + \theta_{i_{k+5}} < \pi$, then
$(i_{k+3}i_{k+4})$ and $(i_{k+4}i_{k+5})$  intersects  
(see figure \ref{Fig:CuspDeform}).  
Let us use the notation  $(i_1i_2i_3)i_4i_5i_6$  and 
$i_1i_2i_3(i_4i_5i_6)$  to indicate 
edges appeared by these perturbations respectively.  

We shall calculate the dihedral angles between this newly intersecting
faces.  
Note that dihedral angles around other (old) edges are
$\pi/2$ by the same argument as in Lemma~\ref{Lem:Intersection}. 

\begin{LemNum}
\label{Lem:DihedralAngle}
Suppose that $\theta=(\theta_1,\ldots,\theta_6)$ be an element of
$\Theta_6$ and $\theta_{i_1}+\theta_{i_2}+\theta_{i_3}<\pi$.  Let
$\omega$ be the dihedral angle between the faces $(i_1 i_2)i_3 i_4 i_5
i_6$ and $i_1(i_2 i_3)i_4 i_5 i_6$. Then we have
$$ \cos\omega = \sqrt{\frac{\sin\theta_{i_1}\sin\theta_{i_3}}
   {\sin(\theta_{i_1}+\theta_{i_2})\sin(\theta_{i_2}+\theta_{i_3})}} $$
\end{LemNum}

To prove this lemma, we shall use the next identity.

\begin{LemNum}
Suppose that $\alpha_1+\alpha_2+\alpha_3+\alpha_4+\alpha_5+\alpha_6
=2\pi$. Then
\begin{equation} \label{eq:hex}
\sin(\alpha_1+\alpha_2)\sin\alpha_4
 - \sin(\alpha_5+\alpha_6)\sin\alpha_3 
 = \sin(\alpha_1+\alpha_2+\alpha_3)\sin(\alpha_3+\alpha_4).
\end{equation}
\end{LemNum}

\begin{pf}
By the addition rule, we have
\begin{align*}
&\sin(\alpha_1+\alpha_2)\sin\alpha_4 - \sin(\alpha_5+\alpha_6)\sin\alpha_3 \\
&=  \sin(\alpha_1+\alpha_2+\alpha_3-\alpha_3)\sin\alpha_4
  - \sin(\alpha_4+\alpha_5+\alpha_6-\alpha_4)\sin\alpha_3 \\
&=  \sin(\alpha_1+\alpha_2+\alpha_3)\cos\alpha_3\sin\alpha_4
   -\cos(\alpha_1+\alpha_2+\alpha_3)\sin\alpha_3\sin\alpha_4 \\
&\phantom{=}
   -\sin(\alpha_4+\alpha_5+\alpha_6)\cos\alpha_4\sin\alpha_3
   +\cos(\alpha_4+\alpha_5+\alpha_6)\sin\alpha_4\sin\alpha_3, \\
\intertext{and by $\sum\alpha_i = 2\pi$,}
&=  \sin(\alpha_1+\alpha_2+\alpha_3)\cos\alpha_3\sin\alpha_4
   -\cos(\alpha_1+\alpha_2+\alpha_3)\sin\alpha_3\sin\alpha_4 \\
&\phantom{=}
   +\sin(\alpha_1+\alpha_2+\alpha_3)\cos\alpha_4\sin\alpha_3
   +\cos(\alpha_1+\alpha_2+\alpha_3)\sin\alpha_4\sin\alpha_3 \\
&= \sin(\alpha_1+\alpha_2+\alpha_3)\sin(\alpha_3+\alpha_4),
\end{align*}
which proves the identity.
\end{pf}

\begin{pf*}{\it Proof of Lemma~\ref{Lem:DihedralAngle}.}
We find the defining equations of two faces as hyperplanes 
in ${\cal E}_{p,\theta}$ in terms of the coordinate $XYZW$.  
We again adopt the notation $s_j, s_{jk},
s_{jkl}$ for $\sin\theta_{i_j}$, $\sin\theta_{i_ji_k}$,
$\sin\theta_{i_ji_ki_l}$.  We first recall the transformation from
${\cal E}_{p,\theta}$ to $XYZW$. We have
\begin{align*}
 Y = \sqrt{\Area B} & = \sqrt{\frac{s_1 s_2}{2s_{12}}} x_{i_1i_2}, \quad
 W = \sqrt{\Area D}   = \sqrt{\frac{s_3 s_4}{2s_{34}}} x_{i_3i_4}, \\
 X = \sqrt{\Area A}
 & = \sqrt{\frac{s_{12} s_{34}}{2 s_{56}}}
     \left(x_{i_2i_3}
           +\frac{s_1}{s_{12}}x_{i_1i_2}
           +\frac{s_4}{s_{34}}x_{i_3i_4}\right) \\
 & = \sqrt{\frac{s_{12} s_{34}}{2 s_{56}}}
     \left(x_{i_2i_3}
           + \sqrt{\frac{2 s_1}{s_2 s_{12}}} \sqrt{\Area B}
           + \sqrt{\frac{2 s_4}{s_3 s_{34}}} \sqrt{\Area D}\right).
\end{align*}
Thus the defining equation of $\{x_{i_1i_2}=0\}$ and
$\{x_{i_2i_3}=0\}$ is 
$$
Y = 0,  \qquad\text{and}\qquad
X =   \sqrt{\frac{s_{34}s_1}{s_{56}s_2}}Y
    + \sqrt{\frac{s_{12}s_4}{s_{56}s_3}}Z.
$$
Hence we can choose normal vectors $n_1, n_2$ for each hyperplane by    
$n_1 = (0,1,0,0)$ and 
$n_2 = (1,\sqrt{\frac{s_{34}s_1}{s_{56}s_2}}, 0, \sqrt{\frac{s_{12}s_4}
{s_{56}s_3}})$ respectively.  Then by
Lemma~\ref{Lem:Distance}
\begin{align*}
\cos\omega
&=\frac{q(n_1, n_2)}{\sqrt{q(n_1, n_1)q(n_2, n_2)}}
 =\frac{\sqrt{\frac{s_{34}s_1}{s_{56}s_2}}}
       {\sqrt{-1+\frac{s_{34}s_1}{s_{56}s_2}+\frac{s_{12}s_4}{s_{56}s_3}}} \\
&=\sqrt{-\frac{s_{56}s_2}{s_{34}s_1}+1+\frac{s_{12}s_4s_2}{s_3s_{34}s_1}}^{-1}
 =\sqrt{1+\frac{s_2}{s_1s_3}\times\frac{s_{12}s_4-s_{56}s_3}{s_{34}}}^{-1}, \\
\intertext{by the identities (\ref{eq:hex}) and (\ref{eq:sin}),}
&=\sqrt{1+\frac{s_2s_{123}}{s_1s_3}}^{-1}
 =\sqrt{1+\frac{s_{12}s_{23}-s_1s_3}{s_1s_3}}^{-1}
 =\sqrt{\frac{s_1s_3}{s_{12}s_{23}}}
\end{align*}
which concludes the proof.
\end{pf*}

As mentioned in the sentences before the
lemma, each dihedral angle around the old edges is $\pi/2$ so they fit
together without producing any singularity.
But the hyperbolic structure at the new edges can be singular  
in  $\overline{X_{6,\theta}}$.  
A cross section perpendicular to the new edge
will be a cone, obtained by taking a 2-dimensional
hyperbolic sector of some angle and identifying the two bounding
rays emanating from the center. 
Such a singular structure appears in the hyperbolic Dehn 
filling theory in  \cite{ThurstonNote} and 
is called a cone singularity.  

\begin{figure}[ht]
\begin{center}
\setlength{\unitlength}{0.00083300in}%
\begingroup\makeatletter\ifx\SetFigFont\undefined%
\gdef\SetFigFont#1#2#3#4#5{%
  \reset@font\fontsize{#1}{#2pt}%
  \fontfamily{#3}\fontseries{#4}\fontshape{#5}%
  \selectfont}%
\fi\endgroup%
\begin{picture}(6900,2385)(226,-2113)
\put(5251,-1861){\makebox(0,0)[lb]{\smash{\SetFigFont{12}{14.4}{\rmdefault}{\mddefault}{\updefault}(56)}}}
\thicklines
\put(1426,-436){\line( 0,-1){375}}
\put(2251,-586){\line( 1, 0){675}}
\put(2551,-436){\line( 0,-1){375}}
\put(3376,-586){\line( 1, 0){675}}
\put(3676,-436){\line( 0,-1){375}}
\put(4501,-586){\line( 1, 0){675}}
\put(4801,-436){\line( 0,-1){375}}
\put(5626,-586){\line( 1, 0){675}}
\put(5926,-436){\line( 0,-1){375}}
\put(5626,-1186){\line( 1, 0){675}}
\put(5926,-1036){\line( 0,-1){375}}
\put(1126,-1186){\line( 1, 0){675}}
\put(1426,-1036){\line( 0,-1){375}}
\put(2251,-1186){\line( 1, 0){675}}
\put(2551,-1036){\line( 0,-1){375}}
\put(3376,-1186){\line( 1, 0){675}}
\put(3676,-1036){\line( 0,-1){375}}
\put(4501,-1186){\line( 1, 0){675}}
\put(4801,-1036){\line( 0,-1){375}}
\put(1426,-211){\line( 0, 1){225}}
\put(2551,-211){\line( 0, 1){225}}
\put(3676,-211){\line( 0, 1){225}}
\put(4801,-211){\line( 0, 1){225}}
\put(5926,-211){\line( 0, 1){225}}
\put(7051,-211){\line( 0, 1){225}}
\put(376,-211){\line( 0, 1){225}}
\put(376,-436){\line( 0,-1){375}}
\put(376,-1036){\line( 0,-1){375}}
\put(376,-1636){\line( 0,-1){150}}
\put(376,-1786){\line( 1, 0){300}}
\put(376,-1186){\line( 1, 0){300}}
\put(376,-586){\line( 1, 0){300}}
\put(7051,-436){\line( 0,-1){375}}
\put(7051,-1036){\line( 0,-1){375}}
\put(7051,-1636){\line( 0,-1){150}}
\put(7051,-1786){\line(-1, 0){300}}
\put(6751,-1186){\line( 1, 0){300}}
\put(6751,-586){\line( 1, 0){300}}
\put(1126,-1786){\line( 1, 0){675}}
\put(1426,-1636){\line( 0,-1){150}}
\put(2251,-1786){\line( 1, 0){675}}
\put(2551,-1636){\line( 0,-1){150}}
\put(3376,-1786){\line( 1, 0){675}}
\put(3676,-1636){\line( 0,-1){150}}
\put(4501,-1786){\line( 1, 0){675}}
\put(4801,-1636){\line( 0,-1){150}}
\put(5626,-1786){\line( 1, 0){675}}
\put(5926,-1636){\line( 0,-1){150}}
\put(1126, 14){\line( 1, 0){675}}
\put(2251, 14){\line( 1, 0){675}}
\put(3376, 14){\line( 1, 0){675}}
\put(4501, 14){\line( 1, 0){675}}
\put(5626, 14){\line( 1, 0){675}}
\put(6751, 14){\line( 1, 0){300}}
\put(376, 14){\line( 1, 0){300}}
\put(601,-1561){564123}
\put(1726,-1561){564213}
\put(2851,-1561){564231}
\put(3976,-1561){564321}
\put(5101,-1561){564312}
\put(6226,-1561){564132}
\put(1276,-1561){(12)}
\put(2401,-1561){(13)}
\put(3526,-1561){(23)}
\put(4651,-1561){(12)}
\put(5776,-1561){(13)}
\put(751,-1861){(56)}
\put(601,-961){546123}
\put(1726,-961){546213}
\put(2851,-961){546231}
\put(3976,-961){546321}
\put(5101,-961){546312}
\put(6226,-961){546132}
\put(1276,-961){(12)}
\put(2401,-961){(13)}
\put(3526,-961){(23)}
\put(4651,-961){(12)}
\put(5776,-961){(13)}
\put(6376,-1261){(46)}
\put(5251,-1261){(46)}
\put(4126,-1261){(46)}
\put(3001,-1261){(46)}
\put(1876,-1261){(46)}
\put(751,-1261){(46)}
\put(601,-361){456123}
\put(1726,-361){456213}
\put(2851,-361){456231}
\put(6226,-361){456132}
\put(1276,-361){(12)}
\put(2401,-361){(13)}
\put(3526,-361){(23)}
\put(5776,-361){(13)}
\put(6376,-661){(45)}
\put(5251,-661){(45)}
\put(4126,-661){(45)}
\put(3001,-661){(45)}
\put(1876,-661){(45)}
\put(226,-361){(23)}
\put(226,-961){(23)}
\put(226,-1561){(23)}
\put(901,-2086){d}
\put(2026,-2086){e}
\put(3151,-2086){f}
\put(4276,-2086){a}
\put(5401,-2086){b}
\put(6526,-2086){c}
\put(826,164){a}
\put(1951,164){b}
\put(3076,164){c}
\put(4201,164){d}
\put(5326,164){e}
\put(6451,164){f}
\put(6901,-1561){(23)}
\put(6901,-961){(23)}
\put(4651,-361){(12)}
\put(5101,-361){456312}
\put(3976,-361){456321}
\put(6901,-361){(23)}
\put(226,-136){x}
\put(226,-736){y}
\put(226,-1336){z}
\put(7126,-1336){z}
\put(7126,-736){y}
\put(7126,-136){x}
\put(751,-661){(45)}
\put(1876,-1861){(56)}
\put(3001,-1861){(56)}
\put(4126,-1861){(56)}
\put(6376,-1861){(56)}
\put(6376,-61){(56)}
\put(5251,-61){(56)}
\put(4126,-61){(56)}
\put(3001,-61){(56)}
\put(1876,-61){(56)}
\put(751,-61){(56)}
\put(1126,-586){\line( 1, 0){675}}
\end{picture}
\end{center}
\caption{Cusp of $\overline{X_6}$: 
	edges with the same letters are identified}\label{Fig:Cusp}
\end{figure}
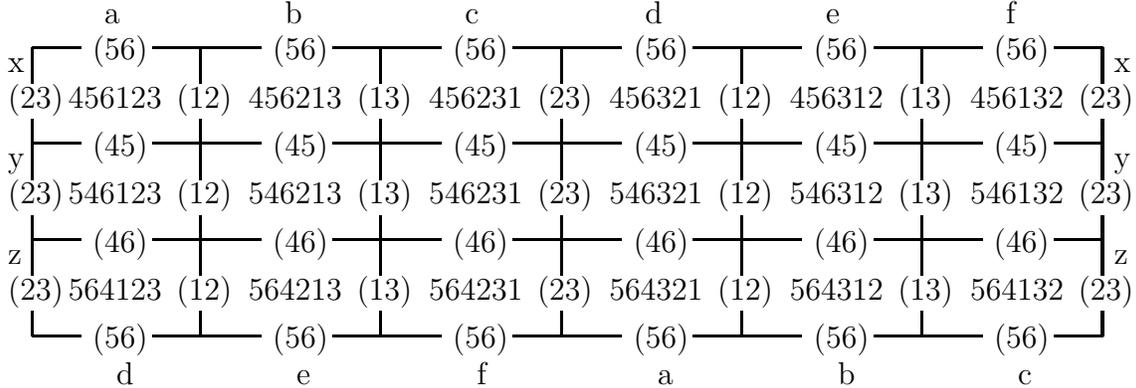

Before investigating what the singularity looks like, 
we describe a polygonal decomposition of cusps.  
The faces of hexahedra  $\Delta_{p,\theta_0}$'s 
lie in the same component of 
cusps in  $\overline{X_6}$  if and only if the labels 
are identical as a partition of six numbers .  
Hence the number of cusps is equal to the number of 
partitions of  $\{ 1, 2, \cdots, 6 \}$  into a pair of 
three numbers,  $= \binom{6}{3}/2 = 10$.  
We may use the notation  $(i_1i_2i_3)(i_4i_5i_6)$  to 
indicate a component of cusps also if there is no confusion.  
Now, since the total amount of ideal vertices of 
hexahedra in  $X_6$  is  
$3 \times (6-1)!/2 = 180$,  
$180/10 = 18$  components of ideal vertices of  
$\Delta_{p,\theta_0}$'s come to a component of 
cusps in  $\overline{X_6}$.  
Figure~\ref{Fig:Cusp}  shows how they come to.  
Replacing the mark  $j$  by  $i_j$  in  Figure~\ref{Fig:Cusp}, 
we get the picture of the cusp labeled by 
$(i_1i_2i_3)(i_4i_5i_6)$.  

To calculate the cone angle for each cone singularity, we
look at the boundary of a small equidistant neighborhood of the
singular locus (or the ideal vertex). 
Note that the boundary is a torus obtained 
by gluing rectangles, each of which 
is a new face of some  $\Delta_{p,\theta}$  appeared by 
the truncation.
Since the old faces of  $\Delta_{p,\theta}$  intersect in right angle, 
not only a combinatorial but conformal 
pattern of its polygonal decomposition is described 
by  Figure~\ref{Fig:Cusp}.  
We then see that six edges with the
same label $(i_5i_6)$ form a nontrivial loop on the boundary. 
If  $\theta_{i_1}+\theta_{i_2}+\theta_{i_3} < \pi$, 
it winds once around the singular locus labeled by  
$(i_1i_2i_3)i_4i_5i_6$.  
If  $\theta_{i_1}+\theta_{i_2}+\theta_{i_3} > \pi$, 
it is homotopic to a loop winding the singular locus 
labeled by  $i_1i_2i_3(i_4i_5i_6)$  twice.  

\begin{LemNum}\label{Lem:Angle}
Suppose that $\{i_1,\ldots,i_6\}=\{1,\ldots,6\}$, $\theta =
(\theta_1, \theta_2, \theta_3, \theta_4, \theta_5, \theta_6) \in
\Theta_6$ and  $\theta_{i_1}+\theta_{i_2}+\theta_{i_3}<\pi$.
Let $\angle((i_1i_2i_3);\theta)$ be the cone
angle about a singular locus labeled by  
$(i_1 i_2 i_3)i_4 i_5i_6$  in  $\overline{X_{6,\theta}}$. 
Then we have
$$ 
\cos(\angle((i_1i_2i_3);\theta)/2)
   = N(i_1i_2i_3;\theta). 
$$
\end{LemNum}

\begin{pf}
Regarding the mark  $j$  as  $i_j$ in Figure~\ref{Fig:Cusp}, 
we can identify the angle with the sum of   
dihedral angles about the new edges of 
six successive hexahedra in the horizontal direction.  
However because of the gluing rule shown in  Figure~\ref{Fig:Cusp},  
the half can be computed by summing three successive ones.  
To compute it, look at 
the faces of hexahedra which appear as sections 
to the singular locus.  
The corresponding Figure~\ref{Fig:Develope-hex} in this case  
to Figure~\ref{Fig:Develope} 
shows how hexahedra are developed about the singular 
locus in  ${\cal P}_6$.  
Choose normal vectors of length $= \sqrt{-1}$ to each face sharing 
the singular locus.  
\begin{figure}[ht]
\begin{center}
\setlength{\unitlength}{0.00083300in}%
\begingroup\makeatletter\ifx\SetFigFont\undefined%
\gdef\SetFigFont#1#2#3#4#5{%
  \reset@font\fontsize{#1}{#2pt}%
  \fontfamily{#3}\fontseries{#4}\fontshape{#5}%
  \selectfont}%
\fi\endgroup%
\begin{picture}(3924,2199)(589,-2548)
\put(2476,-2161){\makebox(0,0)[lb]{\smash{\SetFigFont{12}{14.4}{\rmdefault}{\mddefault}{\updefault}2(13)456}}}
\thicklines
\put(2701,-961){\line( 1, 2){300}}
\put(901,-886){\line( 4, 1){2100}}
\put(3001,-361){\line( 1,-1){1500}}
\put(3751,-736){\vector(-1, 0){825}}
\put(2701,-2236){\line( 0,-1){300}}
\put(1201,-2461){\line( 1, 1){300}}
\put(3901,-2161){\line( 1,-1){300}}
\put(1501,-1261){\line( 4, 1){1200}}
\put(2701,-961){\line( 1,-1){900}}
\put(901,-1411){\line(-4,-1){300}}
\multiput(1501,-1861)(31.91489,31.91489){48}{\makebox(6.6667,10.0000){\SetFigFont{10}{12}{\rmdefault}{\mddefault}{\updefault}.}}
\multiput(3001,-361)(0.00000,-45.53571){29}{\makebox(6.6667,10.0000){\SetFigFont{10}{12}{\rmdefault}{\mddefault}{\updefault}.}}
\put(3826,-811){\makebox(0,0)[lb]{\smash{\SetFigFont{12}{14.4}{\rmdefault}{\mddefault}{\updefault}(123)456 = singular locus}}}
\put(2026,-1861){\makebox(0,0)[lb]{\smash{\SetFigFont{12}{14.4}{\rmdefault}{\mddefault}{\updefault}213456}}}
\put(1201,-2086){\makebox(0,0)[lb]{\smash{\SetFigFont{12}{14.4}{\rmdefault}{\mddefault}{\updefault}(12)3456}}}
\put(751,-1336){\makebox(0,0)[lb]{\smash{\SetFigFont{12}{14.4}{\rmdefault}{\mddefault}{\updefault}1(23)456}}}
\put(1051,-1711){\makebox(0,0)[lb]{\smash{\SetFigFont{12}{14.4}{\rmdefault}{\mddefault}{\updefault}123456}}}
\put(2851,-1861){\makebox(0,0)[lb]{\smash{\SetFigFont{12}{14.4}{\rmdefault}{\mddefault}{\updefault}231456}}}
\put(3526,-2086){\makebox(0,0)[lb]{\smash{\SetFigFont{12}{14.4}{\rmdefault}{\mddefault}{\updefault}(23)1456}}}
\put(1801,-1861){\line( 1, 1){900}}
\put(2701,-961){\line( 0,-1){975}}
\end{picture}
\end{center}
\caption{}\label{Fig:Develope-hex}
\end{figure}
By replacing point vectors on the hyperbola 
in the proof of Lemma~\ref{Lem:CurveLength} by 
normal vectors to the faces in the Minkowski space, 
we can proceed the computation of sum of 
three dihedral angles 
in the light of similarity between Lemma~\ref{Lem:EdgeLength}  
and  Lemma~\ref{Lem:DihedralAngle}.   
The conclusion follows from Lemma~\ref{Lem:Distance} (2a) here 
instead of (1).  
\end{pf}

We very briefly recall some foundations of the 
hyperbolic Dehn filling theory based on Thurston \cite{ThurstonNote}, 
Neumann-Zagier \cite{NeumannZagier} and  
Culler-Shalen \cite{CullerShalen}. 
Let  $N$  be an orientable complete hyperbolic 3-manifold of 
finite volume with  $s$  cusps,  and 
$\rho_0:\pi_1(N)\rightarrow
\SL(2,\complex)$  a lift of the holonomy representation of  $N$. 
The algebro geometric quotient of all $\SL(2, \complex)$-representations 
of  $\Pi = \pi_1(N)$  is called 
a character variety and denoted by  $X(\Pi)$.  
The set of Dehn filled deformations of  $N$  is 
locally parameterized by a neighborhood of 
$\rho_0$  in  $X(\Pi)$.  
The following structure theorem is fundamental, 
which appeared in this form for example in \cite{HodgsonKerckhoff}, 
though the claim could be derived from the arguments in 
\cite{ThurstonNote, CullerShalen, NeumannZagier}.  

\begin{LemNum}\label{Lem:Coordinate}
$X(\pi_1(N))$ is a smooth manifold 
near $\rho_0$ of complex dimension  $s$. 
If $m_1,\ldots,m_s$ are meridional curves for cusps,
then the map $f:X(\pi_1(N))\rightarrow \complex^s$ defined by
$$
f(\chi) = (\chi_\rho(m_1),\ldots,\chi_\rho(m_s)) 
$$  
is a local diffeomorphism near $\rho$ where 
$\chi_\rho(m_i)=\trace \rho(m_i)$.
\end{LemNum}

Going back to our setting and 
let ${\cal L}$ be the union of singular loci (or ideal vertices). 
$\overline{X_{6,\theta}}-{\cal L}$  is 
homeomorphic to  $\overline{X_6}$  and 
$\overline{X_{6,\theta}}$  is its Dehn filled resultant.  
$\overline{X_{6,\theta}}-{\cal L}$ carries a nonsingular but incomplete
hyperbolic metric.  
Let  $\rho_{\theta} : \Pi = \pi_1( \overline{X_{6,\theta}} -
{\cal L}) \rightarrow \SL(2,\complex)$  
be a lift of the holonomy representation 
of  $\overline{X_{6,\theta}} - {\cal L}$.  

To see how  $\overline{X_{6,\theta}}$  is deformed, 
we only need by  Lemma~\ref{Lem:Coordinate}  
to compute the trace of a holonomy 
image of  $\rho_{\theta}$  at some meridional elements.  
To define appropriate meridional elements, 
assume for the moment that  $i_1 + i_2 + i_3 < \pi$.  
Then the cusp labeled by  $(i_1i_2i_3)(i_4i_5i_6)$  becomes  
a singular locus in  $\overline{X_{6,\theta}}$  labeled 
by  $(i_1i_2i_3)i_4i_5i_6$,  
and there is a natural meridional element 
winding once around the singular locus.  
We denote it by  $m_{i_1i_2i_3}$.  
It is homotopic to a loop on the boundary of a tubular 
neighborhood of the singular locus labeled by  
either  $(i_4i_5)$, $(i_5i_6)$  or  $(i_4i_6)$  
(see Figure~\ref{Fig:Cusp}).  
Note that  $m_{i_1i_2i_3}$  is a meridional element if  
$i_1 + i_2 + i_3 < \pi$  but no longer meridional 
if  $i_1 + i_2 + i_3 > \pi$  in any sense.  

For example, $\rho_{\theta}(m_{123})$  in Figure~\ref{Fig:Cusp}  
acts as a translation of six blocks in horizontal direction, 
however, it is actually a rotation, a parabolic translation or 
a hyperbolic translation according to whether  
$\theta_1 + \theta_2 + \theta_3$  is less than, equal to or 
greater than  $\pi$.  
Having this picture in mind, we prove 

\begin{LemNum}\label{Lem:TraceIsNijk}
Suppose that $\{i_1,\ldots,i_6\}=\{1,\ldots,6\}$, $\theta = (\theta_i,
\ldots, \theta_6) \in \Theta_6$, and 
denote by $\rho_\theta$  
the holonomy representation with respect to $\theta$. 
Then we have
$$ \chi_{\rho_\theta}(m_{i_1i_2i_3})=2 N(i_1i_2i_3;\theta). $$
\end{LemNum}

\begin{pf}
Suppose that  $\theta_{i_1} + \theta_{i_2} +
\theta_{i_3} < \pi$  and   
let $\omega$ be $\exp(\sqrt{-1}\angle((i_1i_2i_3)/2;\theta))$.  
Then  $\rho_{\theta}(m_{i_1i_2i_3})$  is an elliptic element 
of rotation $\omega^2$,  and 
its action on $\complex\cup\{\infty\}$ is conjugate to
$$ z \mapsto \omega^2 z
   = \frac{\omega z+0}{0\cdot z+ \omega^{-1}}.
$$
Hence $\rho_\theta(m_{i_1i_2i_3})$ is conjugate to
$$ \begin{pmatrix}
	\omega & 0 \\
	0 & \omega^{-1}
   \end{pmatrix}. 
$$
By Lemma~\ref{Lem:Angle}, 
$$
\chi_{\rho_\theta}(m_{i_1i_2i_3}) = \omega +
\omega^{-1} = 2 \cos (\angle((i_1i_2i_3);\theta)/2) 
= 2 N(i_1i_2i_3;\theta).
$$ 

Suppose that  $\theta_{i_1} + \theta_{i_2} + \theta_{i_3} = \pi$.  
Then  $(i_1i_2i_3)(i_4i_5i_6)$  is an 
ideal vertex  
and  $\rho_{\theta}(m_{i_1i_2i_3})$  is a parabolic element.  
By easy computation, we have 
$$ 
\chi_{\rho_\theta}(m_{i_1i_2i_3}) = 2 = 2 N(i_1i_2i_3;\theta). 
$$

Suppose that  $\theta_{i_1} + \theta_{i_2} + \theta_{i_3} > \pi$.  
Then  $\rho_{\theta}(m_{i_1i_2i_3})$  is a hyperbolic element 
translating the face labeled by, say,  $(i_2i_3)$  to one 
to be identified in  $\overline{X_{6, \theta}}$. 
Regarding  $i_j$  as  $j$  in Figure~\ref{Fig:Cusp}, 
we can identify the faces in question with ones 
on the right hand and left hand sides.  
Let  $\delta$  be a distance between them.  

$\delta$  is equal to the sum of lengths of edges labeled by  
$i_1i_2i_3(i_4i_5i_6)$  in six successive hexahedra in 
the horizontal direction in Figure~\ref{Fig:Cusp}.  
Again because of the gluing rule in  Figure~\ref{Fig:Cusp},  
the half of $\delta$  
can be computed by summing three successive ones. 
Then developing the faces of three hexahedra involving 
the edges labeled by  $i_1i_2i_3(i_4i_5i_6)$  as in 
Figure~\ref{Fig:Develope},  we now realize that 
the computation we carried out in Lemma~\ref{Lem:Angle} 
measures the distance between faces labeled by  
$(i_2,i_3)$  on a side and  $(i_2i_3)$  on the middle 
in Figure~\ref{Fig:Cusp}  because of Lemma~\ref{Lem:Distance} (2b).  
More precisely, we have 
$$
\cosh(\delta/2) = N(i_1i_2i_3; \theta).  
$$
The difference is that the value of Lorentz bilinear forms 
for normal vectors is greater than  $1$  in this case, 
but less than  $1$  in the previous case.  

The action of such a hyperbolic motion is conjugate to 
$$ 
z \mapsto \lambda^2 z
   = \frac{\lambda z+0}{0\cdot z+ \lambda^{-1}},  
$$
where  $\lambda$  is a real number  $> 1$  and 
$\delta = \log \lambda^2$.  
Hence  $\rho_\theta(m_{i_1i_2i_3})$ is conjugate to
$$ \begin{pmatrix}
	\lambda & 0 \\
	0 & \lambda^{-1}
   \end{pmatrix},
$$
and 
$$
\chi_{\rho_\theta}(m_{i_1i_2i_3}) = \lambda + \lambda^{-1} 
= e^{\delta/2} + e^{-\delta/2} 
= 2 \cosh (\delta/2) 
= 2 N(i_1i_2i_3;\theta).
$$ 
This completes the proof.  
\end{pf}

\begin{ThNum}\label{Th:DefromationOfHexagons}
The map from $\Theta_6$ to $X(\Pi)$ defined by 
the assignment $: \theta \mapsto \overline{X_{6,\theta}}$  is 
a local embedding at   
$\theta_0 = (\pi/3, \cdots, \pi/3)$.  
\end{ThNum}
 
\begin{pf}
To apply Lemma~\ref{Lem:Coordinate}, we choose meridional loops for
five cusps, say $m_{145}$, $m_{234}$, $m_{235}$, $m_{245}$, $m_{345}$.
Denote by $\rho_\theta$  a lift of the holonomy of  
$\overline{X_{6,\theta}} - {\cal L}$  in  
$\SL(2,\complex)$.  
Let $\Phi: \Theta_6 \to \complex^5$ be a map defined by
$$ \theta \mapsto 
(\chi_{\rho_\theta}(m_{145}),
 \chi_{\rho_\theta}(m_{234}),
 \chi_{\rho_\theta}(m_{235}),
 \chi_{\rho_\theta}(m_{245}),
 \chi_{\rho_\theta}(m_{345})).
$$
It suffices to show that $\Phi$ is locally injective around $\theta_0$.
Since the image of $\Phi$ is contained in the real part $\real ^5
\subset \complex^5$, we regard $\Phi$ as a map from $\Theta$ to $\real ^5$.
We shall show that the Jacobian of $\Phi$ does not vanish. 

As a basis of the tangent space of $\Theta$ at $\theta_0$, we take the
following five paths.
$p_j(t) = (\theta_{j1}(t), \cdots, \theta_{j6}(t))$ ($j=1, \cdots, 5$)
passing through $\theta_0$ defined by
$$ \theta_{jj}(t)=\dfrac{\pi}{3} -t,\qquad
   \theta_{j6}(t)=\dfrac{\pi}{3} +t,\qquad
   \theta_{jk}(t)=\dfrac{\pi}{3} \quad \text{for}\,  k \ne j, 6.$$

By Lemma~\ref{Lem:TraceIsNijk} we have
$$
\begin{array}{clcccccr}
\Phi(p_1(t)) & = (&f(t), &2,     &2,     &2,     &2     &), \\
\Phi(p_2(t)) & = (&2,    &f(t),  &f(t),  &f(t),  &2     &), \\
\Phi(p_3(t)) & = (&2,    &f(t),  &f(t),  &2,     &f(t)  &), \\
\Phi(p_4(t)) & = (&f(t), &f(t),  &2    , &f(t),  &f(t)  &), \\
\Phi(p_5(t)) & = (&f(t), &2,     &f(t),  &f(t),  &f(t)  &).
\end{array}
$$
where 
\begin{align*}
f(t) 
&= 2 N\left(123;\left(\dfrac{\pi}{3} -t,\dfrac{\pi}{3},\dfrac{\pi}{3},
            \dfrac{\pi}{3},\dfrac{\pi}{3}\right)\right) \\
&= 2\dfrac{\sin\dfrac{\pi}{3}\sin\left(\dfrac{\pi}{3}-t\right)
- \sin\dfrac{\pi}{3}\sin t -2\sin t\sin\left(\dfrac{\pi}{3}-t\right)}%
{\sin^2 \left(\dfrac{2\pi}{3}-t\right)}.
\end{align*}
Then $f'(0) = -6\sqrt{3}.$
Therefore the Jacobian of $\Phi$ at $\theta_0$ is
$-3\times(-6\sqrt{3})^5$ which concludes the proof of the theorem.
\end{pf}

\begin{Remark} 
{\em 
Our computation is valid only when $\theta$  is close to  $\theta_0$.  
\cite{MorinNishi} presents an expanded background for the present deformation
which might help to prove global injectivity.
} 
\end{Remark}

\end{document}